\documentclass{article}

\usepackage[english]{babel}

\usepackage[papersize={8.5in,11in}, margin=1in,footskip=0.4in,vmarginratio={1:1},hmarginratio={1:1}]{geometry}

\usepackage{amsfonts} %extra math fonts
\usepackage{amsmath} %miscellaneous enhancements for math formulas
\usepackage{amssymb} %extra math symbols
\usepackage{amsthm} %theorems set-up
\usepackage{latexsym} %extra math symbols
\usepackage{mathrsfs} %\math­scr command
\usepackage{mathtools}
\usepackage{xcolor}

\numberwithin{equation}{section}
\newtheorem{theorem}{Theorem}[section]

\newtheorem{lemma}[theorem]{Lemma}
\newtheorem{Proposition}[theorem]{Proposition}

\DeclareMathOperator{\curl}{curl}

\newcounter{1}
\newcounter{2}
\newcounter{3}
\newcounter{4}
\title{On the free-boundary Incompressible Elastodynamics with and without surface tension}
\date{}
\author{Longhui Xu}

\begin{document}
	\maketitle
	\begin{abstract}
	We consider a free-boundary problem for the incompressible elastodynamics describing the motion of an elastic medium in a periodic domain with a moving boundary and a fixed bottom under the influence of surface tension. The local well-posedness in Lagrangian coordinates is proved by extending Gu-Luo-Zhang\cite{Gu2024} on incompressible magnetohydrodynamics. We adapt the idea in Luo-Zhang\cite{luo2023compressible} on compressible gravity-capillary water waves to obtain an energy estimate in graphic coordinates. The energy estimate is uniform in surface tension coefficient if the Rayleigh-Taylor sign condition holds and thus yields the zero-surface-tension limit. 	
	\end{abstract}
\tableofcontents
\pagebreak
%-------------------------------------------------------------------------------

	\section{Introduction}
	We consider the incompressible free-boundary elastodynamic equations
	
	\begin{equation} \label{inisys}
		\begin{cases}
			D_tu +\nabla p=\text{div}(\mathbf{F}\mathbf{F}^T)  \qquad & \text{in } \mathcal{D},\\
			\text{div }u=0 \qquad & \text{in } \mathcal{D},\\
			D_t \mathbf{F}=\nabla u \mathbf{F} \qquad & \text{in } \mathcal{D},\\
			\text{div }\mathbf{F}^T=0 \qquad & \text{in } \mathcal{D}.
		\end{cases}
	\end{equation}
	Here $ \mathcal{D}:= \bigcup _{0\leq t\leq T} \{t\} \times \mathcal{D}_t $ and $\mathcal{D}_t:=\{(x_1,x_2)\in \mathbb{T}^2, -b<x_3<\psi(t,x_1,x_2)\}$ is the periodic domain occupied by the elastic medium at time $t$. $u,p$ denotes the velocity and fluid pressure of the elastic medium. $\mathbf{F}:=(\mathbf{F}_{ij})_{3\times3}$ is the deformation tensor and $\mathbf{F}^T:=(\mathbf{F}_{ji})$ is the transpose of $F$. $\mathbf{F}\mathbf{F}^T$ is the Cauchy-Green tensor for neo-Hooken elsatic materials. $\nabla := (\partial_1,\partial_2,\partial_3)$ is the standard spatial derivative, div $X:=\nabla \cdot X= \partial_i X^i$ is the standard divergence for any vector filed $X$. $D_t:= \partial_t +v\cdot \nabla$ is the material derivative. We consider the boundary conditions of (\ref{inisys})
	
	\begin{equation}
		\begin{cases}
			D_t|_{\partial \mathcal {D}}\in\mathcal{T}(\partial \mathcal{D}), \\
			\mathbf{F}^T n =0 \qquad & \text{on } \partial \mathcal{D}_t, \\
			p=\sigma \mathcal {H}  \qquad & \text{on } \partial\mathcal{D}_t.
		\end{cases}
	\end{equation}
	Here $\mathcal{T}(\partial \mathcal {D}_t)$ is the tangent bundle of $\partial \mathcal {D}_t$. The moving boundary $\partial \mathcal {D}_t$ consists of the moving top $\Sigma_t:=\{(x_1,x_2,x_3)\in \mathbb{R}^3: x_3=\psi(t,x_1,x_2)\}$ and the fixed flat bottom $\Sigma_b:=\{(x_1,x_2,x_3)\in \mathbb{R}^3: x_3=-b\}$. The first condition is the kinematic boundary condition indicating the boundary moves with the velocity of the fluid. The second condition shows the deformation tensor is tangential on the boundary, where $n$ is the unit normal vector to the boundary. The third condition shows that the pressure is balanced by surface tension on the moving top $\Sigma_t$, where $\sigma>0$ is the surface tension constant and $\mathcal{H}$ is the mean curvature of the moving surface.

	Now we write the system in terms of columns of the deformation tensor. Let $A_j:=(A_{1j},A_{2j},A_{3j})^T$ be the $j$-th column of $A$ for any $3\times 3$ matrix A. Then the last equation of ($\ref{inisys}$) reads
	\begin{equation}
		0= (\text{div } \mathbf{F}^T)_i = \partial_j(\mathbf{F}^T)_{ij}=\partial_j \mathbf{F}_{ji}=\text{div } \mathbf{F}_i.
	\end{equation}
	It follows that 
	\begin{equation}
		(\text{div } (\mathbf{F}\mathbf{F}^T))_i =\partial_j (\mathbf{F}\mathbf{F}^T)_{ij}=\partial_j(\mathbf{F}_{ik}\mathbf{F}_{jk})=\partial_j \mathbf{F}_{ik} \mathbf{F}_{kj} + \mathbf{F}_{ik}\partial_j \mathbf{F}_{jk}=\partial_j \mathbf{F}_{ik} \mathbf{F}_{kj}=(\mathbf{F}_k\cdot \nabla) \mathbf{F}_k.
	\end{equation}
	Hence, the first equation reads 
	\begin{equation}
		D_tu +\nabla p=(\mathbf{F}_k\cdot \nabla) \mathbf{F}_k.
	\end{equation}
	Since $(\nabla u\mathbf{F}_j)_{i}=\partial_k u_i \mathbf F_{kj}=( \mathbf F_j\cdot \nabla) u_i$, the third equation reads 
	\begin{equation}
		D_t\mathbf{F}_j=( \mathbf F_j\cdot \nabla) u.
	\end{equation}
	Columnwisely , the boundary condition $\mathbf{F}^T n =0$ reads 
	\begin{equation}
		0=(\mathbf{F}^T n)_j=\mathbf{F}_{kj}n_k=\mathbf{F}_j\cdot n.
	\end{equation}
	The system (\ref{inisys}) can be reformulated as 
	\begin{equation}\label{expsys}
		\begin{cases}
			D_tu +\nabla p= (\mathbf{F}_k\cdot \nabla) \mathbf{F}_k \qquad & \text{in } \mathcal{D},\\
			\text{div }u=0 \qquad & \text{in } \mathcal{D},\\
			D_t\mathbf{F}_j=( \mathbf F_j\cdot \nabla) u \qquad & \text{in } \mathcal{D},\\
			\text{div }\mathbf{F}_j=0 \qquad & \text{in } \mathcal{D},
		\end{cases} 
	\end{equation}
	with the boundary condition 
	\begin{equation}\label{expbc}
		\begin{cases}
			D_t|_{\partial \mathcal {D}}\in\mathcal{T}(\partial \mathcal{D}), \\
			\mathbf{F}_j\cdot n =0 \qquad & \text{on } \partial \mathcal{D}_t, \\
			p=\sigma \mathcal {H}  \qquad & \text{on } \partial\mathcal{D}_t.\\
		\end{cases} 
	\end{equation}
	%-------------------------------------------------------------------------------
	
	\subsection{History and Background}
	Let us first briefly review the results of the free-boundary Euler equations, which have been intensively studied in recent decades. The first breakthrough was by Wu \cite{wu1997well,wu1999well} on the local well-posedness(LWP) of the irrotational water wave. For the general incompressible problem with non-zero vorticity, Christodoulou-Lindblad \cite{christodoulou2000motion} established an apriori energy estimate and Lindblad \cite{lindblad2001well,lindblad2005well} proved the LWP for the case without surface tension by Nash-Moser iteration, which leads to a loss of regularity. Coutand-Shkoller \cite{coutand2007well,coutand2010simple} proved the LWP for the case with surface tension by introducing tangential smoothing and avoided the loss of regularity. 
	\par The study of free-boundary compressible fluid is much less and most of the results only deal with the case without surface tension. Lindblad \cite{lindblad2003compwell,lindblad2005compwell} first proved the LWP by Nash-Moser and Trakhinin \cite{trakhinin2009local} extended the LWP to non-isentropic fluid in an unbounded domain, but both results have a loss of regularity. Luo \cite{luo2018motion} established an apriori estimate without loss of regularity for the isentropic case and Luo-Zhang \cite{luo2022local} proved the LWP without using the Nash-Moser iteration. For the case with non-zero surface tension, Luo-Zhang \cite{luo2023compressible} proved the LWP of the gravity-capillary water waves by tangential smoothing and artificial viscosity with an energy estimate uniform in Mach number and the surface tension coefficient, by which the incompressible limit and the zero surface tension limit can be simultaneously justified.
	\par Now we review the developments of incompressible elastodynamics. The fixed-boundary problem of the elastodynamics is well-understood and the global existence is expected and we refer to  \cite{lei2005global,lei2008global,lin2005hydrodynamics,lin2008initial,lin2012some, coutand2006interaction,dafermos2005hyperbolic}. However, the free-boundary problem is more difficult as we have limited boundary regularity and the boundary terms can enter the highest order in the energy estimate. Most of the existing literature neglected the effect of surface tension. Gu-Wang \cite{gu2020well} and Li-Wang-Zhang \cite{li2021well} proved the LWP under a mixed stability condition. Hu-Huang \cite{hu2019well} proved the LWP under the Rayleigh-Taylor sign condition. For the case with surface tension, Gu-Lei \cite{gu2020local} proved the local well-posedness of the two-dimensional case by studying the viscoelastic system in Lagrangian coordinates and the vanishing viscosity limit. 
	\par Studying the case with nonzero surface tension with Lagrangian coordinates is painful due to its complicated description of the boundary. In this paper, we try to adapt the idea in Luo-Zhang \cite{luo2023compressible} to establish an energy estimate for the incompressible elastodynamics with surface tension in graphical coordinates, by which the mean curvature of the boundary can be neatly formulated. Moreover, the energy estimate is uniform in surface tension coefficient and directly leads to zero surface tension limit. However, we have to point out that the LWP in graphical coordinates cannot be directly studied with tangential smoothing, since the divergence-free condition and the boundary condition of the deformation tensor can not be propagated from the constraints on the initial data in the tangentially smoothed system.
	\subsection{Outline of the paper}
	This paper is organised as follows. In Section 2, we first reformulate the system in the Lagrangian coordinates and establish the LWP by extending Gu-Luo-Zhang \cite{Gu2024} on the LWP of Magnetohydrodynamics(MHD), then we introduce the formulation in the graphical coordinates and our main results. In Section 3, we recall some preliminary results that will be useful in the paper. Section 4 is the derivation of our energy estimate and Section 5 concerns the zero surface tension limit.

	%-------------------------------------------------------------------------------
	\section{Reformulation of the system}
	\subsection{Lagrangian Coordinates}
	The system (\ref{expsys})-(\ref{expbc}) is similar to the free-boundary incompressible magnetohydrodynamics considered in Gu-Luo-Zhang\cite{Gu2024} and thus the local wellposedness can be directly proven by extending their result.
	\subsubsection{Reformulation in Lagrangian Coordinates}
	We reformulate the equations in Lagrangian coordinates to transform the free-boundary problem to be a fixed boundary problem on
	\begin{align*}
		\Omega:=\mathbb{T}^2 \times (-b,0).
	\end{align*}
	Let $\Sigma := \mathbb{T}^2\times\{0\}$, $\Sigma_b:=\mathbb{T}^2\times\{-b\}$. and  $\eta:[0,T]\times \Omega \to \mathcal{D}$ be the flow map, i.e.,
	\begin{align*}
		\partial_t\eta(t,x) = u(t,\eta(t,x)), \qquad \eta_0=\eta(0,\cdot),
	\end{align*}
	with
	\begin{align*}
	\eta_0\big|_\Sigma = \psi(0,\cdot)
	\end{align*}

	We introduce the Lagrangian variables 
	\begin{align*}
	v(t,x)=u(t,\eta(t,x)),\qquad q(t,x)=p(t,\eta(t,x)),\qquad F(t,x)=\mathbf{F}(t,\eta(t,x)).
	\end{align*}
	We define the cofactor matrix $a=[\nabla \eta]^{-1}$, it follows that the material derivative $D_t$ reduces to time derivative $\partial_t$ and the spacial derivative $\partial_j$ becomes $\nabla_a^j:=a^{ij}\partial_i$. Let $N$ be the unit outer normal vector of $\partial \Omega$, the system (\ref{expsys})-(\ref{expbc}) becomes the following:
	
	\begin{equation}\label{Lagsys}
		\begin{cases}
			\partial_t v +\nabla_a q= (F_k\cdot \nabla_a) F_k \qquad & \text{in } \Omega,\\
			\nabla_a \cdot v=0 \qquad & \text{in } \Omega,\\
			\partial_t F_j=(  F_j\cdot \nabla_a) v \qquad & \text{in } \Omega,\\
			\nabla_a\cdot F_j=0 \qquad & \text{in } \Omega
		\end{cases} 
	\end{equation}
	with the boundary condition 
	\begin{equation}\label{Lagbc}
		\begin{cases}
			v_3=F_{3j}=0 \qquad & \text{on } \Sigma_b\\
			a^{\mu \nu}F_{j\nu }N_\mu=0  \qquad & \text{on } \Sigma, \\
			a^{\mu\nu}N_\mu q=-\sigma (\sqrt{g}\Delta_g\eta^\nu)  \qquad & \text{on } \Sigma.\\
		\end{cases}
	\end{equation}
	where $g$ is the metric induced on $\Sigma_t=\eta(t,\Sigma)$ by the embedding $\eta $ and $\Delta_g$ is the Laplacian of $g$. 	
	\subsubsection{Local Well-posedness in Lagrangian Coordinates}
	 This system (\ref{Lagsys})- (\ref{Lagbc}) is identical to the one considered in Gu-Luo-Zhang\cite{Gu2024} by replacing $F_j$ by the magnetic field $b$, hence we can apply their result to obtain the local well-posedness, provided the initial data satisfies certain compatibility conditions. 
	 \begin{theorem}
	 	(Local existence) Let $v_0,\eta_0\in H^{4.5}(\Omega)\cap H^5 (\Sigma) $ and $F_j^0\in H^{4.5}$ be divergence-free vector fields with $(F_j^0\cdot N)|_\Sigma =0$ and define the initial data $q_0$ of $q$ to satisfies the following elliptic equation
	 	\begin{equation}
	 		\begin{cases}
	 		-\Delta q_0=(\partial v_0)(\partial v_0)	-(\partial F_j^0)(\partial F_j^0)\qquad & \text{in } \Omega,\\
	 		q_0=\sigma \mathcal{H}_0 \qquad & \text{on } \Sigma,\\
	 		\frac{\partial q_0}{\partial N}=0 \qquad & \text{on } \Sigma_b.
	 		\end{cases}
	 	\end{equation}
	 	Then there exists some $T>0$ depending on $\sigma, v_0, F_j^0$, such that the system (\ref{expsys})-(\ref{expbc}) with initial data $(v_0,F_j^0,q_0)$ has a unique strong solution $(\eta, v,q)$ with the energy estimates
	 	\begin{align}
	 	\sup_{0\leq t\leq T}E(t)\leq C,
	 	\end{align}
	 	where $C$ depends on $||v_0||_{4.5},||F^0_J||_{4.5},|v_0|_5$ and
	 	\begin{align}
	 	E(t):=&||\eta(t)||^2_{4.5}+\sum _{k=0}^3\bigg(||\partial_t^kv(t)||^2_{4.5-k},	||\partial_t^k(F_j^0\cdot \partial)\eta(t)||^2_{4.5-k}\bigg)+||\partial^4_tv(t)||_0^2+||\partial^4_t(F_j^0\cdot \partial)\eta(t)||_0^2\\
	 	&+\sum_{k=0}^3\big(\big|\bar{\partial}(\Pi\bar{\partial}^k\partial_t^{3-k}v(t)\big)\big|_0^2\big)+\big|\bar{\partial}(\Pi\bar{\partial}^3 (F_j^0\cdot \partial)\eta(t)\big)\big|_0^2,
	 	\end{align}
	 	here $\Pi$ is the canonical normal projection defined on the moving interface.
	 \end{theorem}
%-------------------------------------------------------------------------------

	\subsection{Graphic Coordinates}
	The remaining parts of this paper are devoted to developing an $\sigma-$uniform energy estimate, provided the Rayleigh-Taylor sign condition is satisfied, and thus we can obtain the zero-surface tension limit. This allow us to show the LWP for the case of $\sigma = 0$ 
	\subsubsection{Reformulation in Graphic Coordinates}
	We introduce graphical coordinates to convert the free-boundary problem into a fixed domain problem on 
	\begin{equation*}
		\Omega:=\mathbb{T}^2 \times (-b,0).
	\end{equation*}
	Let $\overline{x}=(x_1,x_2)$. The moving top $\Sigma_t:=\{(\overline{x},x_3)\in \Omega: x_3=\psi(t,\overline{x})\}$ is represented by the graph of $\psi$. To fix the interior, we introduce the extension of $\psi$
	\begin{equation}
		\varphi(t,\overline{x},x_3)=x_3+\chi(x_3)\psi(t,\overline{x}),
	\end{equation}
	where $\chi\in C_c^\infty(-b,0]$ satisfying 
	\begin{equation}\label{chi}
		||\chi'||_{L^\infty (-b,0]}\leq \frac{1}{1+||\psi_0||_\infty}, \qquad \chi\equiv 1 \text{ on }(-\delta_0,0]
	\end{equation}
	for some small constant $\delta_0>0$.\\
	Note that
	\begin{equation}\label{p3phi}
		\partial_3 \varphi=1+\chi'(x_3)\psi(t,\overline{x}).
	\end{equation}
	By (\ref{chi}), we have 
	\begin{equation}\label{p3phi2}
		|\chi'(x_3)\psi(t,\overline{x})|\leq\frac{||\psi(t)||_\infty}{1+||\psi_0||_\infty}<1, \qquad t\in[0,T]
	\end{equation}
	for some small $T>0$. It follows from (\ref{p3phi}) and (\ref{p3phi2})
	\begin{equation}
		\partial_3 \varphi (t,\overline{x},x_3)\geq c_0 >0,  \qquad t\in[0,T] \label{partial3phi}
	\end{equation} 
	for some $c_0>0$. Let
	\begin{equation}
		\Phi(t,\overline{x},x_3) = (\overline{x},\varphi(t,\overline{x},x_3)),
	\end{equation}
	then (\ref{partial3phi}) ensures us to define the diffeomorphism $\displaystyle \Phi(t,\cdot) :\Omega \to \mathcal{D}_t$.
	
	Now we introduce the graphical variables
	\begin{equation}
		v(t,x)=u(t,\Phi(t,x)), \qquad q(t,x)=p(t,\Phi(t,x)), \qquad F(t,x)=\mathbf{F}(t,\Phi(t,x)).
	\end{equation}
	Also, we introduce the induced differential operator $\displaystyle \partial ^\phi$ so that for any function $g(t,x), x\in\mathcal{D}_t$ and $G:=g(t,\Phi(t,x))$, we have
	\begin{equation}
		\partial_\alpha^\varphi G=\partial_\alpha g (t,\Phi(t,x)),  \qquad \alpha=t,1,2,3.
	\end{equation}
	It follows that 
	\begin{align}
		\partial_t^\varphi&=\partial_t-\frac{\partial_t \varphi}{\partial_3 \varphi}\partial_3,\\
		\nabla_\tau^\varphi=\partial_\tau^\varphi&=\partial_\tau-\frac{\partial_\tau \varphi}{\partial_3 \varphi}\partial_3, \qquad \tau=1,2,\\
		\nabla_3^\varphi= \partial_3^\varphi&=\frac{1}{\partial_3 \varphi}\partial_3.
	\end{align}
	We define the cofactor matrix
	\begin{equation}
		A:=
		\begin{pmatrix}
			1 & 0 & 0 \\
			0 & 1 & 0 \\
			-\frac{\partial_1 \varphi}{\partial_3 \varphi} 
			& -\frac{\partial_2 \varphi}{\partial_3 \varphi} & \frac{1}{\partial_3 \varphi}
		\end{pmatrix}.
	\end{equation}
	It follows that 
	\begin{equation}
		\nabla_i^\varphi=A_{ji} \partial_j.
	\end{equation}
	Next, we introduce the material derivative 
	\begin{align}
		D_t^\varphi&:=\partial_t^\varphi+v\cdot\nabla^\varphi \\
		&= (\partial_t-\frac{\partial_t \varphi}{\partial_3 \varphi}\partial_3)
		+v_1(\partial_1-\frac{\partial_1 \varphi}{\partial_3 \varphi}\partial_3)
		+v_2(\partial_1-\frac{\partial_2 \varphi}{\partial_3 \varphi}\partial_3)
		+\frac{v_3}{\partial_3 \varphi}\partial_3 \nonumber\\
		&=\partial_t +\overline{v}\cdot \overline{\partial}
		+(-\partial_t \varphi-v_1\partial_1\varphi 
		-v_2\partial_2\varphi+v_3)\frac{1}{\partial_3\varphi}\partial_3
		\nonumber\\
		&=\partial_t+\overline{v}\cdot\overline{\partial}+(v\cdot\mathbf{N} -\partial_t\varphi)\partial_3^\varphi \label{refmart},
	\end{align}
	where $\mathbf{N}:=(-\partial_1\varphi,-\partial_2\varphi,1)$ and $\overline{\partial}=\overline{\nabla}=(\partial_1,\partial_2)$. 

	Now, we reformulate the boundary conditions. Let $N:=(1,0,\partial_1 \psi)\times(0,1,\partial_2 \psi) = (-\partial_1 \psi,-\partial_2\psi,1)$ be the standard normal vector to the moving top $\Sigma_t$ and $n=(0,0,1)$ the unit normal vector to the fixed bottom. By (\ref{refmart}),
	\begin{equation}
		D_t^\varphi|_{\partial\Omega}= \partial_t+\overline{v}\cdot\overline{\partial}+\big((v\cdot\mathbf{N} -\partial_t\varphi)\partial_3^\varphi\big)|_{\Sigma\cup\Sigma_b}.
	\end{equation}
	since 
	\begin{equation}
		(v\cdot\mathbf{N} -\partial_t\varphi)|_{\Sigma}=v\cdot N -\partial_t \psi,
	\end{equation}
	and
	\begin{equation}
		(v\cdot\mathbf{N} -\partial_t\varphi)|_{\Sigma_b}=v\cdot n-\partial_t(-b)=v_3,
	\end{equation}
	we have $D_t^\varphi$ is tangential if
	\begin{align}
		v\cdot N -\partial_t \psi&=0 \qquad \text{ on } \Sigma,\\
		v_3&=0 \qquad\text{ on } \Sigma_b.
	\end{align}
	Next, for the boundary condition of $q$, by definition of mean curvature,
	\begin{align}
		q=-\sigma \overline{\nabla}\cdot\big(\frac{\overline{\nabla\psi}}{\sqrt{1+|\overline{\nabla}\psi|^2}}) &\qquad\text{ on } \Sigma.
	\end{align}
	Consequently the systems (\ref{inisys}) and (\ref{expsys}) are converted into 
	\begin{equation}\label{phisys}
		\begin{cases}
			D_t^\varphi v +\nabla^\varphi q= (F_k\cdot \nabla^\varphi) F_k \qquad & \text{in } \Omega,\\
			\nabla^\varphi \cdot v=0 \qquad & \text{in } \Omega,\\
			D_t^\varphi F_j=(F_j\cdot \nabla^\varphi) v  \qquad & \text{in } \Omega,\\
			\nabla^\varphi \cdot  F_j=0 \qquad & \text{in } \Omega
		\end{cases} 
	\end{equation}
	with the boundary conditions
	\begin{equation}\label{phibc}
		\begin{cases}
			\partial_t \psi=v\cdot N &\qquad \text{ on } \Sigma,\\
			q=-\sigma \overline{\nabla}\cdot\big(\frac{\overline{\nabla\psi}}{\sqrt{1+|\overline{\nabla}\psi|^2}}) &\qquad\text{ on } \Sigma,\\
			F_j\cdot N=0 &\qquad \text{ on } \Sigma,\\
			v_3=0 &\qquad\text{ on } \Sigma_b,\\
			F_{3j}=0 &\qquad \text{ on } \Sigma_b,\\
			q=0 &\qquad\text{ on } \Sigma_b.\\
		\end{cases}
	\end{equation}
	We can show that div $F_j=0$ and the boundary condition $F_j\cdot N=0$ are only constraints on the initial data.
	\begin{Proposition} 
	Let the initial data satisfies
	\begin{align} \label{inicons1}
	\nabla^\varphi  \cdot F_j =0	,
	\end{align}
	and the boundary condition
	\begin{align}\label{inicons2}
	F_j\cdot n=0 \quad\text{on } \Sigma\cup\Sigma_b.
	\end{align}
	Then the solution satisfies (\ref{inicons1}) and (\ref{inicons2}) for all $t\in[0,T].$
	\end{Proposition}
	\noindent\textbf{Proof}
	Taking $\nabla^\varphi \cdot$ on $D_t^\varphi F_j-(F_j\cdot \nabla^\varphi) v=0$, we have 
	\begin{align*}
		0&=\partial_i^\varphi(\partial_t^\varphi F_{ij}+\big(v\cdot \nabla^\varphi)F_{ij}\big)-\partial_i^\varphi\big((F_k \cdot \nabla^\varphi )v_i\big)\\
		&=D_t^\varphi \big(\nabla^\varphi \cdot F_j\big)+(\partial_i^\varphi v\cdot \nabla^\varphi)F_{ij}-(\partial_i^\varphi F_k\cdot \nabla^\varphi )v_i
		-(F_k\cdot \nabla^\varphi )(\nabla^\varphi \cdot v),
	\end{align*}
	where the fourth term vanishes due to the divergence-free condition of $v$ and the third term cancels the second term. We have 
	\begin{align}
		D_t^\varphi \big(\nabla^\varphi \cdot F_j\big) =0,
	\end{align}
	which is a linear equation of $\nabla^\varphi \cdot F_j$ and the propagation of $\nabla^\varphi \cdot F_j$ follows from standard characteristic curve method. For the boundary condition, we consider $D_t^\varphi F_j-(F_j\cdot \nabla^\varphi) v+F_j(\nabla^\varphi \cdot v)=0$ on the boundary, let $\tau=1,2$
	\begin{align*}
	0&=(\partial_t +\overline{v}\cdot \overline{\partial}) F_j-F_{\tau j}(\partial_\tau-\partial_\tau\varphi\partial_3)v-F_{3j}\partial_3 v+F_j(\partial_\tau-\partial_\tau\varphi\partial_3)v_\tau+F_j\partial_3 v_3\\
	&=(\partial_t +\overline{v}\cdot \overline{\partial}) F_j-F_{\tau j}\partial_\tau v-(F_j\cdot N)\partial_3 v+F_j\partial_\tau v_\tau+F_j(N\cdot \partial_3 v).
	\end{align*}
	We times $N$ on both sides to have 
	\begin{align*}
		0&=(\partial_t +\overline{v}\cdot \overline{\partial}) F_j\cdot N-F_{\tau j}(\partial_\tau v\cdot N)-(F_j\cdot N)(\partial_3 v\cdot N)+(F_j\cdot N)\partial_\tau v_\tau+(F_j\cdot N)(N\cdot \partial_3 v)\\
		&=(\partial_t +\overline{v}\cdot \overline{\partial})F_j\cdot N-F_{\tau j}(\partial_\tau v\cdot N)+(F_j\cdot N)\partial_\tau v_\tau\\
		&=(\partial_t +\overline{v}\cdot \overline{\partial})(F_j\cdot N)-F_j\cdot \partial_t N-F_j\cdot (\overline{v}\cdot \overline{\partial})N-F_{\tau j}(\partial_\tau v\cdot N)+(F_j\cdot N)\partial_\tau v_\tau,
	\end{align*}
	where the second term cancels exactly the third and fourth term,
	\begin{align*}
	-F_j\cdot \partial_t N=F_{\tau j}	 \partial_t \partial_\tau \psi=F_{\tau j}(\partial_\tau v\cdot N+v\cdot \partial_\tau N)=F_{\tau j}(\partial_\tau v\cdot N)+F_{\tau j}(\overline{v}\cdot \overline{\partial}) N.
	\end{align*}
	Hence, we have a linear equation of $F_j\cdot N$ on the boundary
	\begin{align}
	0=	(\partial_t +\overline{v}\cdot \overline{\partial})(F_j\cdot N)+\partial_\tau v_\tau(F_j\cdot N).
	\end{align}

	\subsection{The main result}
	\begin{theorem}
	Let $\sigma >0$, $(v,q,F,\psi)$ satisfies the above system. Define 
	\begin{align}
	E(t)=\sum_{k=0}^4\bigg(||\partial_t ^kF||_{4-k},||\partial_t ^kv||_{4-k},||\sqrt{\sigma}\partial_t ^k \overline{\nabla}\psi||_{4-k}\bigg)	+ \sum _{k=0}^3||\partial_t^k q||_{4-k},
	\end{align}
	Then there exists $T>0$ such that 
	\begin{align}
	E(t)\leq P(E(0))+\int_0^T	P(E(t)).
	\end{align}
	Furthermore, if the Rayleigh-Taylor sign condition $-\partial_3 q\geq c_0>0$ is assumed,
	 then $\sum _{k=0}^4|\partial_t^k\psi|_{4-k}$ enters the energy and the estimate become $\sigma-$uniform.

	\end{theorem}

	\section{The auxiliary results}
	
	\subsection{Preliminary Lemmas}
	\begin{lemma}
		For $\alpha,\beta =t,1,2,3$ and generic function $F$, 
		\begin{equation}
			[\partial_\alpha^\varphi,\partial_\beta^\varphi]F=0.
		\end{equation}
	\end{lemma}
	\noindent\textbf{Proof}
	It directly follows from the fact that $\partial_\alpha^\varphi F= \partial_\alpha f (t,\Phi(t,x))$, where $F=f(t,\Phi(t,x))$. More precisely
	\begin{align*}
		\partial_\alpha^\varphi(\partial_\beta^\varphi F) = \partial_\alpha^\varphi\big(\partial_\beta f(t,\Phi(t,x))\big)= \partial_\alpha\partial_\beta f(t,\Phi(t,x))=\partial_\beta\partial_\alpha f(t,\Phi(t,x))= \partial_\beta^\varphi(\partial_\alpha^\varphi F).
	\end{align*}
	
	\begin{lemma}\label{ibp}
		(Integration by parts) Let $g=g(t,x), f=f(t,x), x\in \Omega$,
		\begin{align}
		\int_\Omega 	(\partial_i^\varphi f)g\partial_3 \varphi =-\int_\Omega 	 f(\partial_i^\varphi g)\partial_3 \varphi +\int_\Sigma fgN_i +\int_{\Sigma_b} fgn_i.
		\end{align}
	where $N=(-\partial_1\psi,-\partial_2 \psi ,1) ,n=(0,0,1)$ are the normal vector function of $\Sigma,\Sigma_b$ respectively.
	\end{lemma}
		\noindent\textbf{Proof} The proof follows from change of variable and standard integration by parts. Let $G(t,x), F(t,x) ,x\in \mathcal{D}_t $ such that $g(t,x) =G(t,\Phi(t,x)) ,f=F(t,\Phi(t,x))$, then
	\begin{align}
	\int_\Omega 	(\partial_i^\varphi f)g\partial_3 \varphi &=\int_\Omega \partial_i F(t,\Phi(t,x))G(t,\Phi(t,x))\partial_3 \varphi \nonumber\\
	&=\int _{\mathcal{D}_t}   \partial_iF(t,x)G(t,x)\nonumber\\
	&=-\int_{\mathcal{D}_t}F(t,x)\partial_iG(t,x) +\int _{\partial\mathcal{D}_t}F(t,x)G(t,x) \mathbf{n}_i \nonumber\\
	&=-\int_\Omega 	 f(\partial_i^\varphi g)\partial_3 \varphi +\int_\Sigma fgN_i +\int_{\Sigma_b} fgn_i \nonumber.
	\end{align}

	\begin{lemma}\label{transp}
		(Transport Theorem) Let $f=f(t,x), x\in \Omega$
		\begin{align}
		\frac{d}{dt}\int_\Omega f\partial_3 \varphi =\int_\Omega D_t^\varphi f\partial_3 \varphi .
		\end{align}
	\end{lemma}
	\noindent\textbf{Proof}
	The proof follows from definition of $D_t^\varphi$, $\nabla^\varphi \cdot v =0$ and integration by parts,
	\begin{align}
	\frac{d}{dt}\int_\Omega f\partial_3 \varphi 
	&=\int_\Omega \partial_t f\partial_3 \varphi +  \int_\Omega f\partial_t\partial_3 \varphi \nonumber \\
	&=\int_\Omega \partial_t^\varphi f \partial_3 \varphi +\int_\Omega \partial_t \varphi \partial_3 f+  \int_\Omega f\partial_t\partial_3 \varphi \nonumber\\
	&=\int_\Omega D_t^\varphi f \partial_3 \varphi -\int_\Omega (v\cdot \nabla^\varphi ) f \partial_3 \varphi  +\int _\Sigma f\partial_t\varphi \nonumber\\
	&=\int _\Omega D_t^\varphi \partial_3 \varphi  +\int_\Omega (\nabla^\varphi \cdot v)f\partial_3 \varphi \nonumber\\
	&= \int_\Omega D_t^\varphi \partial_3 \varphi  \nonumber.
	\end{align}

	\subsection{Elliptic Estimates}
	
	\begin{lemma}\label{hodge}
		\textbf{(The Hodge-type elliptic estimate) } Let $X$ be a smooth vector field and $s\geq 1$, then 
		\begin{equation}
			||X||_s^2\lesssim C\big(|\psi|_{s},|\overline{\partial}\psi|_{W^{1,\infty}}\big)\big(||\nabla^\varphi\cdot X||_{s-1}^2+||\nabla^\varphi\times X||_{s-1}^2+||\overline{\partial}^s X||_0^2+||X||_0^2\big).\\
		\end{equation}
	\end{lemma}
	We refer to lemma B.2 of \cite{ginsberg2020local} for the proof

	\begin{lemma}\label{ellipest}
		\textbf{(Low regularity elliptic estimate)} 
		Assume $W \in H^1$ with $W=0$ on $\Sigma_b$ satisfying
		\begin{equation}\label{elli}
			\begin{cases}
				-\Delta^\varphi W = \nabla^\varphi \cdot \pi \qquad & \text{in }\Omega,\\
				\nabla ^\varphi W\cdot N =h \qquad & \text{on } \partial \Omega,
			\end{cases}
		\end{equation}
		where $\pi,\nabla^\varphi\cdot\pi \in L^2(\Omega)$ and $h\in H^{-0.5}(\partial \Omega)$ with the compatibility condition 
		\begin{equation}
			\int_{\partial \Omega} (\pi \cdot N+h)dS=0.
		\end{equation}
		Then, we have 
		\begin{equation}
			||W||_1\lesssim _{vol(\Omega)} ||\partial_3 \varphi||_\infty||\pi||_0+|\pi \cdot N+h|_{-0.5}.
		\end{equation}
	\end{lemma}
	\noindent\textbf{Proof}
	Testing the first equation of (\ref{elli}) with $\omega\partial_3 \phi, \omega\in H^1$, 
	\begin{equation}\label{ibp:1}
		\int _\Omega (-\partial_i^\varphi \partial_i ^\varphi W)\omega \partial_3 \phi 
		= \int_\Omega (\partial_i ^\varphi \pi) \omega \partial_3 \phi.
	\end{equation}
	By integration by parts, 
	\begin{align}
		\int _\Omega (-\partial_i^\varphi \partial_i ^\varphi W)\omega \partial_3 \phi 
		\nonumber
		&= \int _\Omega (\partial_i ^\varphi W)\partial_i^\varphi\omega \partial_3 \phi -  \int _{\partial \Omega} (\partial_i ^\varphi W)\omega N_i\\
		\nonumber
		&=\int _\Omega A^{\mu i} A^{\nu i} \partial _\mu W \partial_\nu\omega \partial_3 \phi -  \int _{\partial \Omega} h\omega,\\
		\nonumber
		\int_\Omega (\partial_i ^\varphi \pi) \omega \partial_3 \phi 
		&= -\int_\Omega \pi \cdot \nabla^\varphi \omega \partial_3 \phi 
		+\int_{\partial \Omega} (\pi \cdot N) \omega.
		\nonumber
	\end{align}
	It follows from (\ref{ibp:1}) and $\partial_3 \phi>0$ that 
	\begin{align}
		\int _\Omega A^{\mu i} A^{\nu i} \partial _\mu W \partial_\nu\omega 
		\leq \int _\Omega A^{\mu i} A^{\nu i} \partial _\mu W \partial_\nu\omega \partial_3 \phi
		&= -\int_\Omega \pi \cdot \nabla^\varphi \omega \partial_3 \phi 
		+\int_{\partial \Omega} (\pi \cdot N+h) \omega 
		\nonumber \\
		&\lesssim ||\pi||_0||\partial_3 \phi ||_\infty ||\omega||_1+|\pi\cdot N+h|_{-0.5}|\omega|_{0.5}.
	\end{align}
	Since $\displaystyle AA^T$ is positive definite and $||AA^T||_\infty\leq K$, by elliptic estimate 
	\begin{equation}
		||\partial W||_0\lesssim ||\partial_3 \phi||_\infty||\pi||_0 + |\pi\cdot N+h|_{-0.5}.
	\end{equation}
	By Poincar\'e's inequality
	\begin{equation}
		||W||_0
		\lesssim _{vol(\Omega)} ||\partial W||_0 +\int_\Omega  W,
	\end{equation}
	where 
	\begin{equation}
		\int_\Omega W dx= \int_\Omega \partial_3 x_3 W dx= -\int_\Omega x_3\partial_3 W dx\leq ||x_1||_0||\partial W||_0\lesssim_{Vol(\Omega)}||\partial W||_0,
	\end{equation}
	the boundary terms vanishes since $x_3 = 0 $ on $\Sigma$ and $W=0$ on $\Sigma_b$.
	\section{Energy estimate}
	\subsection{Pressure Estimate}
	By the first equation of (\ref{phisys}),
	\begin{equation}
		-\nabla^\varphi q= D_t^\varphi v - (F_k\cdot \nabla^\varphi) F_k \label{momentum}.
	\end{equation}
	Taking $\displaystyle \nabla^\varphi\cdot $, we have
	\begin{equation}
		-\Delta^\varphi q= \nabla^\varphi \cdot (D_t^\varphi v - (F_k\cdot \nabla^\varphi) F_k)  \qquad \text{ on } \Omega
	\end{equation}
	with the boundary condition 
	\begin{equation}
		\nabla^\varphi q \cdot N= -(D_t^\varphi v - (F_k\cdot \nabla^\varphi) F_k) \cdot N \qquad \text{ on } \Sigma\cup\Sigma_b,
	\end{equation}
	which is exactly the form of the elliptic system in lemma (\ref{ellipest}) and satisfies the compatibility condition. It follows that 
	\begin{equation}
		||q||_1\lesssim_{vol(\Omega)}||\partial_3 \varphi||_\infty (||D_t^\varphi v||_0 + ||F_k||_0||\nabla^\varphi F_k||_0 ).
	\end{equation}
	
	For the higher order term, by definition or $\partial_3^\varphi$ and taking the normal component of (\ref{momentum}),
	\begin{align}
	\partial_3 q =\partial_3 \varphi \partial_3 ^\varphi q= -\partial_3 \varphi D_t^\varphi v +\partial_3 \varphi(F_k\cdot \nabla^\varphi )F_{3k}.
	\end{align}
	Let $\partial_\gamma= \partial_t, \partial_1,\partial_2,\partial_3$ and $k=1,2,3$. We have
	\begin{align}
	||\partial_\gamma^k	 \partial_3 q||_0=||\partial_\gamma^k	\big(-\partial_3 \varphi D_t^\varphi v +\partial_3 \varphi(F_k\cdot\nabla^\varphi )F_{3k}\big)||_0\leq P\bigg(\sum_{l=0}^3 |\partial_t ^l \overline{\nabla}\psi|_{3-l},\sum _{l=0}^4 ||\partial_t^k v||_{4-k}, \sum _{l=0}^4 ||\partial_t^k F||_{4-k}\bigg).
	\end{align}
	Next, by definition of $\partial_\tau^\varphi, \tau=1,2$,
	\begin{align}
	\partial_\tau q	=\partial_\tau^\varphi q+\frac{1}{\partial_3 \varphi}\partial_3 q=D_t^\varphi v_\tau -(F_k\cdot \nabla^\varphi )F_{\tau k} +\frac{1}{\partial_3 \varphi}\partial_3 q.
	\end{align}
	Let $D=\overline{\partial}$ or $\partial_t,k=1,2,3$, we have 
	\begin{align}
		||D^k\partial_\tau q||_0=||D^k\big(D_t^\varphi v_\tau -(F_k\cdot \nabla^\varphi )F_{\tau k} +\frac{1}{\partial_3 \varphi}\partial_3 q\big)||_0\leq P\bigg(\sum_{l=0}^3 |\partial_t ^l \overline{\nabla}\psi|_{3-l},\sum _{l=0}^4 ||\partial_t^k v||_{4-k}, \sum _{l=0}^4 ||\partial_t^k F||_{4-k}\bigg).
	\end{align}

	\subsection{Div-Curl Estimate}
	We adopt Hodge-type elliptic estimates to study 
	\begin{equation*}
		||\partial_t^k v||_{4-k}^2, \qquad ||\partial_t^kF_j||_{4-k}^2 \qquad \text{for } k=0,1,2,3,4.
	\end{equation*}
	By lemma (\ref{hodge}), replacing $\displaystyle X$ by $\partial_t^k v$ and $\displaystyle \partial_t^k F_j$ and $s$ by $3-k$, we have 
	\begin{align}
		||\partial_t^k v||_{4-k}^2&\lesssim C\big(|\psi|_{4-k},|\overline{\partial}\psi|_{W^{1,\infty}}\big)\big(||\nabla^\varphi\cdot \partial_t^k v||_{3-k}^2+||\nabla^\varphi\times \partial_t^k v||_{3-k}^2+||\overline{\partial}^{4-k} \partial_t^k v||_0^2+|| v||_0^2\big),\\
		||\partial_t^kF_j||_{4-k}^2&\lesssim C\big(|\psi|_{4-k},|\overline{\partial}\psi|_{W^{1,\infty}}\big)\big(||\nabla^\varphi\cdot \partial_t^kF_j||_{3-k}^2+||\nabla^\varphi\times \partial_t^kF_j||_{3-k}^2+||\overline{\partial}^{4-k} \partial_t^k F_j||_0^2+||F_j||_0^2\big).
	\end{align}
	
	\subsubsection{$L^2$-estimates}
	Testing $v\partial_3\varphi$ with the first equation of (\ref{phisys}), we have
	\begin{align}\label{l2test}
		\int_\Omega (D_t^\varphi v +\nabla^\varphi q) \cdot v\partial_3 \varphi= \int_\Omega (F_k\cdot \nabla^\varphi) F_k \cdot v\partial_3 \varphi.
	\end{align}
	On the left-hand side, by the transport theorem (\ref{transp}) and integration $\nabla^\varphi$ by parts (\ref{ibp}),
	\begin{align}
		\int_\Omega (D_t^\varphi v +\nabla^\varphi q) \cdot v\partial_3 \varphi = \frac{1}{2}\frac{d}{dt}\int_\Omega |v|^2\partial_3 \varphi -\int_\Omega q (\nabla^\varphi \cdot v) \partial_3 \varphi + \int_\Sigma q (v\cdot N),
	\end{align}
	where the boundary term on $\Sigma_b$ vanishes due to the slip boundary condition $\displaystyle v\cdot n=v_3=0$ on  $\Sigma_b$ and the second term $\displaystyle \int_\Omega q(\nabla^\varphi \cdot v)\partial_3 \varphi$ vanishes due to the incompressible condition $\nabla^\varphi\cdot v=0$. \\
	Plugging in the kinematic boundary condition and the boundary condition for $q$, then integrate $\overline{\nabla}$ by parts
	\begin{align}
		\int_\Sigma q v\cdot N
		&=-\sigma \int_\Sigma \overline{\nabla} \cdot\frac{\overline{\nabla}\psi}{\sqrt{1+|\overline{\nabla} \psi |^2}} \partial_t\psi\\
		&=\sigma\int_\Sigma  \frac{\overline{\nabla}\psi}{\sqrt{1+|\overline{\nabla} \psi |^2}} \cdot \overline{\nabla} \partial_t\psi\\
		&=\sigma \frac{d}{dt}\int_\Sigma \sqrt{1+|\overline{\nabla}\psi|^2}.
	\end{align}\\
	On the right-hand side of (\ref{l2test}), since $\nabla^\varphi \cdot F_j=0$ we can integrate $(F_k\cdot \nabla^\varphi)$ by parts to get 
	\begin{align}\label{RHS:F}
		\int_\Omega (F_k\cdot \nabla^\varphi) F_k \cdot v\partial_3 \varphi 
		&= -\int_\Omega  F_k \cdot (F_k\cdot \nabla^\varphi) v\partial_3 \varphi,
	\end{align}
	where the boundary term vanishes since $F_j\cdot N=0$ on $\Sigma\cup \Sigma_b.$ Since $D_t^\varphi F_j= (F_j \cdot \nabla^\varphi) v$, it follows from (\ref{RHS:F}) that 
	\begin{align}
		\int_\Omega (F_k\cdot \nabla^\varphi) F_k \cdot v\partial_3 \varphi =-\int_\Omega  F_k \cdot D_t^\varphi F_k \partial_3 \varphi =-\frac{1}{2}\frac{d}{dt}\sum_{k=1}^3\int_\Omega |F_k|^2\partial_3 \varphi.
	\end{align}
	Therefore,
	\begin{equation}
		\frac{d}{dt}\bigg(\frac{1}{2}\sum_{k=1}^3\int_\Omega |F_k|^2\partial_3 \varphi +\frac{1}{2}\int_\Omega |v|^2\partial_3 \varphi+\sigma \int_\Sigma \sqrt{1+|\overline{\nabla}\psi|^2}\bigg)=0.
	\end{equation}

	\subsubsection{Curl estimates}
	
	In this section, we control $\displaystyle ||\nabla^\varphi\times \partial_t^k v||_{3-k}^2$ and $||\nabla^\varphi\times \partial_t^kF_j||_{3-k}^2$. Recall in the system (\ref{phisys}) that we have
	\begin{align}
		D_t^\varphi v +\nabla^\varphi q&= (F_k\cdot \nabla^\varphi) F_k, \\
		D_t^\varphi F_j&=(F_j\cdot \nabla^\varphi) v. 
	\end{align}
	Taking the curl operator, we have 
	\begin{align}
		D_t^\varphi (\curl^\varphi v) 
		&= (F_k\cdot \nabla^\varphi) (\curl^\varphi F_k)+J_1 \label{curl:0.1},\\
		D_t^\varphi (\curl^\varphi F_j)&=(F_j\cdot \nabla^\varphi)(\curl^\varphi v) + K_1 \label{curl:0.2},
	\end{align}
	where $J_1:=[\curl^\varphi,D_t^\varphi]v
	+[\curl^\varphi,(F_k\cdot \nabla^\varphi)]F_k$  and $K_1:=[\curl^\varphi, D_t^\varphi]F_j+ [\curl^\varphi, (F_j\cdot \nabla^\varphi)]v$. Next taking $\partial^\alpha$, $\alpha = (\alpha_1, \alpha_2, \alpha_3)$ with $\alpha\leq 3$, 
	\begin{align}\label{equ:2}
		D_t^\varphi (\partial^\alpha \curl^\varphi v) 
		&= (F_k\cdot \nabla^\varphi) (\partial^\alpha \curl^\varphi F_k)+\partial^\alpha J_1 +J_2,\\
		D_t^\varphi (\partial^\alpha\curl^\varphi F_j)&=(F_j\cdot \nabla^\varphi)(\partial^\alpha \curl^\varphi v) + \partial^\alpha K_1+K_2,
		\label{curl:1}
	\end{align}
	where $J_2:=[\partial^\alpha, D^\varphi_t]\curl^\varphi v
	+[\partial^\alpha, (F_k\cdot \nabla^\varphi)]\curl^\varphi F_k$ and 
	$K_2 =[\partial^\alpha,D_t^\varphi]\curl^\varphi F_j+[\partial^\alpha, (F_k\cdot \nabla^\varphi)]\curl^\varphi v$.
	
	Now we test (\ref{equ:2}) with $\partial^\alpha \curl^\varphi v \partial_3 \varphi$
	\begin{align}
		\frac{1}{2}\frac{d}{dt}\int_\Omega|\partial^\alpha \curl^\varphi v|^2\partial_3\varphi
		= \int_\Omega (F_k\cdot \nabla^\varphi) (\partial^\alpha \curl^\varphi F_k)\cdot(\partial^\alpha \curl^\varphi v) \partial_3 \varphi
		+\int_\Omega (\partial^\alpha J_1 +J_2)\cdot(\partial^\alpha \curl^\varphi v) \partial_3 \varphi.
		\label{curl:2}
	\end{align}
	Then we integrate $(F_k\cdot \nabla^\varphi)$ by parts to get 
	\begin{align}
		&\int_\Omega (F_k\cdot \nabla^\varphi) (\partial^\alpha \curl^\varphi F_k)\cdot(\partial^\alpha \curl^\varphi v) \partial_3 \varphi\\
		&=-\int_\Omega  (\partial^\alpha \curl^\varphi F_k)\cdot(F_k\cdot \nabla^\varphi)(\partial^\alpha \curl^\varphi v) \partial_3 \varphi 
		\nonumber-\int_\Omega (\nabla^\varphi\cdot F_k) (\partial^\alpha \curl^\varphi F_k)\cdot(\partial^\alpha \curl^\varphi v) \partial_3 \varphi
		\nonumber\\
		&=-\int_\Omega  (\partial^\alpha \curl^\varphi F_k)\cdot D_t^\varphi (\partial^\alpha\curl^\varphi F_j)\partial_3 \varphi
		\nonumber+\int_\Omega  (\partial^\alpha \curl^\varphi F_k)(\partial^\alpha K_1+K_2 ) \partial_3 \varphi
		\nonumber\\
		&\quad-\int_\Omega (\nabla^\varphi\cdot F_k) (\partial^\alpha \curl^\varphi F_k)\cdot(\partial^\alpha \curl^\varphi v) \partial_3 \varphi
		\nonumber\\
		&=-\frac{1}{2}\frac{d}{dt}\int_\Omega  |\partial^\alpha\curl^\varphi F_k|^2\partial_3 \varphi
		\nonumber+\int_\Omega  (\partial^\alpha \curl^\varphi F_k)(\partial^\alpha K_1+K_2 )\partial_3 \varphi
		\nonumber\\
		&\quad-\int_\Omega (\nabla^\varphi\cdot F_k) (\partial^\alpha \curl^\varphi F_k)\cdot(\partial^\alpha \curl^\varphi v) \partial_3 \varphi,
		\label{curl:3}
	\end{align}
	where the boundary term vanishes as  $F_j\cdot N=0.$ It follows from (\ref{curl:2}) and (\ref{curl:3}) that 
	\begin{align}
		\frac{1}{2}\frac{d}{dt}\int_\Omega(|\partial^\alpha \curl^\varphi v|^2+|\partial^\alpha \curl^\varphi F_k|^2)\partial_3\varphi
		=&\int_\Omega (\partial^\alpha J_1 +J_2)\cdot(\partial^\alpha \curl^\varphi v) \partial_3 \varphi
		\nonumber\\
		&+\int_\Omega  (\partial^\alpha K_1+K_2 )(\partial^\alpha \curl^\varphi F_k)\partial_3 \varphi
		\nonumber\\
		&-\int_\Omega (\nabla^\varphi\cdot F_k) (\partial^\alpha \curl^\varphi F_k)\cdot(\partial^\alpha \curl^\varphi v) \partial_3 \varphi.
		\label{curl}
	\end{align}
	It remains to control the error terms on the right-hand side of (\ref{curl}). We first control $\partial^\alpha J_1$. Since
	\begin{align*}
		\curl^\varphi(D_t^\varphi v)_i &=\epsilon^{i\alpha\beta}\partial_\alpha^\varphi(D_t^\varphi )v_\beta\\
		&=\epsilon^{i\alpha\beta}\partial_\alpha^\varphi(\partial_t^\varphi +v_k\partial_k^\varphi )v_\beta\\
		&=\epsilon^{ij\beta} (\partial_t^\varphi +v_k\partial_k^\varphi )\partial_\alpha^\varphi v_\beta +\epsilon^{i\alpha\beta} \partial_\alpha^\varphi v_k \partial_k^\varphi v_\beta,
	\end{align*}
	we have
	\begin{equation} \label{error:1}
		[\curl^\varphi, D_t^\varphi]v = \epsilon^{i\alpha\beta} \partial_\alpha^\varphi v_k \partial_k^\varphi  v_\beta.
	\end{equation}
	Also
	\begin{align*}
		\curl^\varphi\big( (F_k\cdot \nabla^\varphi)F_k\big)_i 
		&=\epsilon ^{i\alpha\beta}\partial^\varphi_\alpha \big( (F_k\cdot \nabla^\varphi) F_{\beta k}\big)\\
		&=\epsilon^{i\alpha\beta}(F_k \cdot\nabla^\varphi) \partial_\alpha^\varphi F_{\beta k} 
		+\epsilon ^{i\alpha\beta}  \big( (\partial^\varphi_\alpha F_k) \cdot \nabla^\varphi\big ) F_{\beta k},
	\end{align*}
	it follows that 
	\begin{equation}\label{error:2}
		[\curl^\varphi,(F_k \cdot \nabla^\varphi)]F_k =\epsilon ^{i\alpha\beta}  \big( (\partial^\varphi_\alpha F_k) \cdot \nabla^\varphi\big ) F_{\beta k}.
	\end{equation}
	Since the terms of (\ref{error:1}) and (\ref{error:2}) are up to the first order, they can be easily controlled by the energy even after taking $\partial^\alpha$. For $J_2$, since $D_t^\varphi$ only contains the first order of $\varphi$ and $\curl^\varphi$ has only 1 derivative on $v$, $[\partial^\alpha, D_t^\varphi]\curl^\varphi v$ can have at most $4$ derivative on $\varphi$ and $v$. More precisely, up to lower order terms,
	\begin{align*}
		\partial^\alpha(D_t^\varphi \curl^\varphi v) 
		&=\partial^\alpha \big( (\partial_t +\overline{v}\cdot\overline{\partial}+(v\cdot\mathbf{N} -\partial_t \varphi)\frac{1}{\partial_3 \varphi}\partial_3 \big) \curl^\varphi v\\
		&\stackrel{\text{L}}{=}D_t^\varphi(\partial^\alpha \curl^\varphi v) + (v\cdot \partial^\alpha \mathbf{N}-\partial^\alpha\partial_t\varphi) \frac{1}{\partial_3\varphi}\partial_3 \curl^\varphi v -(v\cdot\mathbf{N} -\partial_t \varphi) \frac{\partial^\alpha\partial_3 \varphi}{(\partial_3\varphi)^2}\partial_3 \curl^\varphi v,
	\end{align*}
	where $\partial^\alpha \mathbf{N}=(-\partial^\alpha\partial_1 \varphi,-\partial^\alpha\partial_2 \varphi, 0)$. The reasoning for $[\partial^\alpha, (F_k\cdot \nabla^\varphi)]$ is the same. Hence,
	\begin{equation}
		\int_\Omega (\partial^\alpha J_1 +J_2)\cdot(\partial^\alpha \curl^\varphi v) \partial_3 \varphi
		\lesssim P(||v||_4,||F_j||_4,|\psi|_4,|\partial_t \psi|_3).
	\end{equation}
	The argument for the control of $\int_\Omega  (\partial^\alpha K_1+K_2 )(\partial^\alpha \curl^\varphi F_k)\partial_3 \varphi$ is similar by replacing $v$ by $F_k$.
	
	Now, for $\displaystyle\int_\Omega (\nabla^\varphi\cdot F_k) (\partial^\alpha \curl^\varphi F_k)\cdot(\partial^\alpha \curl^\varphi v) \partial_3 \varphi$, again, we can see that the terms of $(\partial^\alpha \curl^\varphi F_k)$ and $(\partial^\alpha \curl^\varphi v)$ are up to forth order of $F_k,v$ and $\varphi$. We can then control the highest order term by $L^2$, 
	\begin{equation}
		\int_\Omega (\nabla^\varphi\cdot F_k) (\partial^\alpha \curl^\varphi F_k)\cdot(\partial^\alpha \curl^\varphi v) \partial_3 \varphi\lesssim P(||v||_4,||F||_4,|\psi|_4).
	\end{equation}
	Let $k=1,2,|\alpha|=3-k$. The estimate for the general cases for $\displaystyle ||\nabla^\varphi\times \partial_t^k v||_{3-k}^2$ and $||\nabla^\varphi\times \partial_t^kF_j||_{3-k}^2$ follows from a parallel argument by taking $\partial^\alpha\partial_t^k$  on $(\ref{curl:0.1})$ and $(\ref{curl:0.1})$ with the number of time derivatives of $\varphi$ up to 4. Consequently,
	\begin{equation}
		\frac{1}{2}\frac{d}{dt}\int_\Omega(|\partial^\alpha\partial_t^k\curl^\varphi v|^2+|\partial^\alpha\partial_t^k \curl^\varphi F_k|^2)\partial_3\varphi
		\lesssim P(||\partial_t^kv||_{4-k}, ||\partial_t^k F_j||_{4-k}, |\partial_t^k \psi |_{4-k} ).
	\end{equation}

	\subsection{Tangential Estimate}
	To get the tangential estimate $||\overline{\partial}^{4-k}\partial_t^k v||_0^2 $, we have to take $D^4=\overline{\partial}^{4-k}\partial_t^k$ on the first equation of (\ref{phisys}), the process produces a trouble term $D^4 \partial_\tau ^\varphi q, \tau=1,2$. Since $\partial_\tau^\varphi$ contains $\partial_\tau \varphi$, the commutator $[D^4, \partial_\tau  ^\varphi] q$ contains $D^4\partial_\tau \varphi$, a term with 4 derivatives on $\varphi$, which can be controlled by $|D^4\partial_\tau \psi|_0$. However, the highest energy term of $\varphi$ from the tangential esimtae is $|\sqrt{\sigma} D^4 \overline{\nabla}\psi  |_0$. It follows that the commutator $[D^4, \partial_\tau  ^\varphi] q$ cannot be controlled $\sigma-$uniformly. 
	
	\subsubsection{Alinhac's Good Unknown}
	To avoid the appearance of the top order of $\varphi$, we introduce the Alinhac's good unknown. Let $D^\alpha= \partial_t^{\alpha_0}\partial_1^{\alpha_1}\partial_2^{\alpha_2}$, $|\alpha|\leq4$,
	\begin{align}
		D^\alpha \partial_\tau^\varphi f
		&=D^\alpha (\partial_\tau f  -\frac{\partial_\tau \varphi}{\partial_3 \varphi} \partial_3 f)\nonumber \\
		&=\partial^\varphi_\tau D^\alpha f - \frac{\partial_\tau D^\alpha  \varphi}{\partial_3 \varphi}\partial_3 f +\frac{\partial_\tau \varphi \partial_3 D^\alpha \varphi}{(\partial_3 \varphi)^2}\partial_3 f +\mathcal{C}_\tau '(f) \nonumber \\
		&=\partial^\varphi_\tau D^\alpha f - (\partial_\tau -\frac{\partial_\tau \varphi }{\partial_3 \varphi}\partial_3) D^\alpha \varphi\partial_3^\varphi f  +\mathcal{C}_\tau '(f) \nonumber\\
		&=\partial^\varphi_\tau(D^\alpha f-D^\alpha \varphi\partial_3^\varphi f) +\mathcal{C}_\tau(f),
	\end{align}
	where for $\tau=1,2$ and $|\beta|=1$
	\begin{align}
		\mathcal{C}_\tau'(f)&=-[D^\alpha,\frac{\partial_\tau \varphi}{\partial_3 \varphi },\partial_3 f]-\partial_3 f[D^\alpha,\partial_\tau \varphi,\frac{1}{\partial_3 \varphi}]+\partial_3 f \partial _\tau \varphi[D^{\alpha-\beta},\frac{1}{(\partial_3 \varphi)^2}]D^\beta \partial_3 \varphi,\\
		\mathcal{C}_\tau (f)&= D^\alpha \varphi\partial_\tau^\varphi\partial_3^\varphi f +\mathcal{C}_\tau'(f) \label{Ctau}.
	\end{align}
	Then by similar calculation, we get 
	\begin{align}
		D^\alpha \partial_3^\varphi f 
		&=D^\alpha ( \frac{1}{\partial_3 \varphi}\partial_3 f ) \nonumber \\
		&=\partial_3 ^\varphi D^\alpha f -\frac{\partial_3 D^\alpha  \varphi}{(\partial_3 \varphi)^2}\partial_3 f +\mathcal{C}_3'(f) \nonumber\\
		&=\partial_3 ^\varphi D^\alpha f -\partial_3^\varphi D^\alpha  \varphi\partial_3^\varphi f +\mathcal{C}_3'(f) \nonumber\\
		&=\partial_3 ^\varphi (D^\alpha f - D^\alpha  \varphi\partial_3^\varphi f) +D^\alpha \varphi (\partial_3 ^\varphi)^2 f +\mathcal{C}_3'(f)\nonumber\\
		&=\partial_3 ^\varphi (D^\alpha f - D^\alpha  \varphi\partial_3^\varphi f)+\mathcal{C}_3(f) \label{C3},
	\end{align}
	where 
	\begin{align}
		\mathcal{C}_3'(f)&=[D^\alpha,\frac{1}{\partial_3 \varphi},\partial_3 f]- \partial_3 f [D^{\alpha-\beta},\frac{1}{(\partial_3 \varphi)^2}]D^\beta \partial_3 \varphi,\\
		\mathcal{C}_3(f)&=D^\alpha \varphi (\partial_3 ^\varphi)^2 f +\mathcal{C}_3'(f).
	\end{align}
	Next, we study 
	\begin{align}
		D^\alpha D_t^\varphi f
		&= D^\alpha (\partial_t +\overline v\cdot \overline{\nabla} +\frac{1}{\partial_3 \varphi}(v\cdot \mathbf{N}-\partial_t \varphi)\partial_3)f\nonumber\\
		&=D_t^\varphi D^\alpha f -\frac{\partial_3 D^\alpha \varphi}{(\partial_3 \varphi)^2}	(v\cdot \mathbf{N} -\partial_t \varphi)\partial_3 f +\frac{1}{\partial_3 \varphi} (v\cdot D^\alpha\mathbf{N}-\partial_t D^\alpha \varphi )\partial_3 f+\mathcal{D}'(f)\nonumber\\
		&=D_t^\varphi D^\alpha f-\big[\frac{1}{\partial_3 \varphi} (v\cdot \mathbf{N} -\partial_t \varphi)\partial_3 D^\alpha  \varphi   +(\overline{v}\cdot \overline{\nabla}D^\alpha \varphi +\partial_tD^\alpha \varphi)\big]\partial_3^\varphi f+\mathcal{D}'(f)\nonumber\\
		&=D_t^\varphi D^\alpha f-D_t^\varphi D^\alpha \varphi\partial_3 ^\varphi f+\mathcal{D}'(f)\nonumber\\
		&=D_t^\varphi (D^\alpha f- D^\alpha \varphi\partial_3 ^\varphi f)+D^\alpha \varphi D_t^\varphi \partial_3^\varphi f +\mathcal{D}'(f)\nonumber\\
		&=D_t^\varphi (D^\alpha f- D^\alpha \varphi\partial_3 ^\varphi f)+\mathcal{D}(f)\label{D},
	\end{align}
	where
	\begin{align}
		\mathcal{D}'(f)&=[D^\alpha, \overline{v}]\cdot \overline{\nabla}f+[D^\alpha, \frac{1}{\partial_3 \varphi}(v\cdot \mathbf N-\partial_t \varphi),\partial_3 f]+[D^\alpha, \frac{1}{\partial_3 \varphi},v\cdot \mathbf{N}-\partial_t \varphi]\partial_3 f\nonumber\\
			&\quad-(v\cdot \mathbf{N} -\partial_t \varphi)\partial_3 f[D^{\alpha-\beta},\frac{1}{(\partial_3 \varphi)^2}]D^\beta \partial_3 \varphi + \frac{1}{\partial_3 \varphi}\partial_3 f[D^\alpha,v]\mathbf{N},\\
			\mathcal{D}(f)&=D^\alpha \varphi D_t^\varphi \partial_3^\varphi f +\mathcal{D}'(f).
	\end{align}
	The quantity $\displaystyle \mathbf{F}:=D^\alpha f- D^\alpha \varphi\partial_3 ^\varphi f$ is called the Alinhac's Good Unknown of $f$. It follows that the control of $||F||_0$ yields the control for $||D^\alpha f||_0$, more precisely
	\begin{align}
	||D^\alpha f||_0\leq ||F||_0 + ||\partial_3^\varphi f||_0||D^\alpha \varphi||_0\label{AGU}.
	\end{align}

	\subsubsection{Reformulation in Alinhac's Good Unknown}
	
	Let 
	\begin{align}
	\mathbf{V}=D^\alpha v - D^\alpha \varphi\partial_3^\varphi v,\quad
	\mathbf{F}=D^\alpha F - D^\alpha \varphi\partial_3^\varphi F,\quad
	\mathbf{Q}=D^\alpha q - D^\alpha \varphi\partial_3^\varphi q.
	\end{align}
	Taking $D^\alpha$ on the first equation of (\ref{phisys}), then by (\ref{Ctau}), (\ref{C3}) and (\ref{D}),
	\begin{align}
	D_t^\varphi \mathbf{V}_i + \partial_i^\varphi \mathbf{Q}&=D^\alpha \big((F_k\cdot \nabla^\varphi)F_{ik}\big) -\mathcal{D}(v_i)-\mathcal{C}_i(q)\nonumber\\
	&=F_{lk}D^\alpha \partial_l^\varphi F_{ik} +[D^\alpha, F_{lk}] \partial_l^\varphi F_{ik}-\mathcal{D}(v_i)-\mathcal{C}_i(q)\nonumber\\
	&=(F_k\cdot \nabla^\varphi) \mathbf{F}_{ik}+F_{lk}C_l(F_{ik})+[D^\alpha, F_{lk}] \partial_l^\varphi F_{ik}-\mathcal{D}(v_i)-\mathcal{C}_i(q)\nonumber\\
	&=(F_k\cdot \nabla^\varphi) \mathbf{F}_{ik}+\mathcal{R}_i^1\label{A1},
	\end{align}
	where 
	\begin{align}
	\mathcal{R}_i^1=F_{lk}\mathcal{C}_l(F_{ik})+[D^\alpha, F_{lk}] \partial_l^\varphi F_{ik}-\mathcal{D}(v_i)-\mathcal{C}_i(q)\label{R1}.
	\end{align}
	Similarly, we have 
	\begin{align}
	D_t^\varphi \mathbf{F}_{ij}&=(F_j\cdot \nabla^\varphi)\mathbf{V}_i+\mathcal R_{ij}^2 \label{A2},\\
	\partial_k^\varphi \mathbf{V}_k&=-\mathcal{C}_k(v_k) \label{A3},
	\end{align}
	where 
	\begin{align}
	\mathcal{R}_{ij}^2=F_{kj}\mathcal{C}_k(v_i)+[D^\alpha,F_{kj}]\partial_k v_i-\mathcal{D}(F_{ij})	\label{R2}.
	\end{align}
	
	Next, we reformulate the boundary conditions. By definition of $\mathbf{Q}$,
	\begin{align}
	\mathbf{Q}|_\Sigma
	&=(D^\alpha q-D^\alpha \varphi \partial_3^ \varphi q)|_\Sigma \nonumber\\
	&=-\sigma D^\alpha\big(\overline{\nabla}\cdot \frac{\overline{\nabla}\psi}{|N|}\big)-D^\alpha \psi \partial_3 q \label{bcq} \qquad \text{on }\Sigma.
	\end{align}
	For the kinematic boundary condition, on $\Sigma$, 
	\begin{align}
	\partial_t D^\alpha \psi 
	&=D^\alpha(v\cdot N)\nonumber\\
	&=D^\alpha (-\overline{v}\cdot \overline{\partial}\psi +v_3)\nonumber\\
	&=-D^\alpha \overline{v}\cdot \overline{\partial}\psi-\overline{v}\cdot D^\alpha \overline{\partial}\psi +D^\alpha v_3 -\sum_{\substack{|\beta_1|+|\beta_2|=3\\ |\beta_1|,|\beta_2|>0}}	 D^{\beta_1}\overline{v}\cdot D^{\beta_2}\overline{\partial}\psi \nonumber\\
	&=D^\alpha v\cdot N-\overline{v}\cdot D^\alpha \overline{\partial}\psi-\sum_{\substack{|\beta_1|+|\beta_2|=3\\ |\beta_1|,|\beta_2|>0}}	 D^{\beta_1}\overline{v}\cdot D^{\beta_2}\overline{\partial}\psi \nonumber\\
	&=\mathbf{V}\cdot N-\overline{v}\cdot D^\alpha \overline{\partial}\psi+\mathcal{S}_1 \label{kbc},
	\end{align}
	where
	\begin{align}
	\mathcal{S}_1:=D^\alpha \psi \partial _3 v\cdot N -\sum_{\substack{|\beta_1|+|\beta_2|=4\\ |\beta_1|,|\beta_2|>0}}	 D^{\beta_1}\overline{v}\cdot D^{\beta_2}\overline{\partial}\psi \label{S1}.
	\end{align}
	For the slip boundary condition, since $D^\alpha \varphi|_{\Sigma_b}=0$ and $v_3|_{\Sigma_b}=0$ on $\Sigma_b$, 
	\begin{align}
		\mathbf{V}_3=\mathbf{V}\cdot n=D^\alpha v \cdot n=0 \label{sbc}.
	\end{align}

	\subsubsection{Energy estimate with full spatial derivatives}
	In this section, we study the equations with full spacial derivative, more precisely, $D^\alpha=\partial_1^{\alpha_1}\partial_2^{\alpha_2}$. 
	
	\begin{theorem}\label{spatial-estimate}
	There exists $T>0$ so that
	\begin{align}
	||D^\alpha v_i||_0^2	+||D^\alpha F_{ik}||_0^2+|\sqrt{\sigma}D^\alpha \overline{\nabla}\psi|_0^2\lesssim P(E(0))+\int_0^T P(E(t)).
	\end{align}
	Furthermore, if the Rayleigh-Taylor sign condition $\-\partial_3 q\geq c_0 >0$ is assumed on $\Sigma$, the term $|D^\alpha\psi|_0^2$ enters the energy and the estimate become $\sigma$-uniform
	\begin{align}
		||D^\alpha v_i(T)||_0^2	+||D^\alpha F_{ik}||_0^2+|\sqrt{\sigma}D^\alpha \overline{\nabla}\psi|_0^2+|D^\alpha\psi|_0^2\lesssim P(E(0))+\int_0^T P(E(t)).
	\end{align}
	\end{theorem}
	Testing (\ref{A1}) with $\mathbf{V}_i\partial_3 \varphi$, we have 
	
	\begin{align}
	\frac{1}{2}\frac{d}{dt}\int_\Omega |\mathbf{V}_i|^2\partial_3 \varphi 
	&= \int_\Omega \mathbf{V}_iD_t^\varphi\mathbf{V}_i \partial_3 \varphi \nonumber\\
	&=-\int_\Omega\partial_i^\varphi \mathbf{Q} \mathbf{V}_i \partial_3 \varphi + \int_\Omega(F_k\cdot \nabla^\varphi) \mathbf{F}_{ik}\mathbf{V}_i \partial_3 \varphi+\int_\Omega \mathcal{R}_i^1\mathbf{V}_i \partial_3 \varphi \label{interior1}.
	\end{align}
	Integrating $\partial_i^\varphi$ by parts, then by $\mathbf{V}\cdot n=0$ on $\Sigma_b$ and (\ref{A3}),
	\begin{align}
		-\int_\Omega\partial_i^\varphi \mathbf{Q} \mathbf{V}_i \partial_3 \varphi
		&=\int _\Omega \mathbf{Q} \partial_i^\varphi \mathbf{V}_i \partial_3 \varphi -\int_\Sigma \mathbf{Q}\mathbf{V}\cdot N \nonumber\\
		&=-\int_\Omega \mathcal{C}_i(v_i)\mathbf{Q}\partial_3 \varphi-\int_\Sigma \mathbf{Q}\mathbf{V}\cdot N \label{interior2}.
	\end{align}
	Next, integrating $\nabla^\varphi$ by parts, then by (\ref{A2}) and $\nabla^\varphi \cdot F_k =0$,
	\begin{align}
	\int_\Omega(F_k\cdot \nabla^\varphi) \mathbf{F}_{ik}\mathbf{V}_i \partial_3 \varphi 
	&= -\int_\Omega \mathbf{F}_{ik}(F_k\cdot \nabla^\varphi)\mathbf{V}_i \partial_3 \varphi-\int_\Omega (\nabla ^\varphi \cdot F_k)\mathbf{F}_{ik}\mathbf{V}_i \partial_3 \varphi \nonumber\\
	&= -\int_\Omega \mathbf{F}_{ik}D_t^\varphi \mathbf{F}_{ik} \partial_3 \varphi +\int_\Omega \mathbf{F}_{ik} \mathcal{R}^2_{ik} \partial_3 \varphi\nonumber\\
	&= -\frac{1}{2}\frac{d}{dt} \int_\Omega |\mathbf{F}_{ik}|^2 \partial_3 \varphi +\int_\Omega \mathbf{F}_{ik} \mathcal{R}^2_{ik} \partial_3 \varphi\label{interior3},
	\end{align}
	where the boundary term vanishes due to $F_k\cdot N=0$ on $\Sigma \cup \Sigma_b$.
	Summing up (\ref{interior1}), (\ref{interior2}) and (\ref{interior3}), we have 
	\begin{align}
		\frac{1}{2}\frac{d}{dt}(\int_\Omega |\mathbf{V}_i|^2\partial_3 \varphi +\int_\Omega |\mathbf{F}_{ik}|^2 \partial_3 \varphi)=
			-\int_\Sigma \mathbf{Q}\mathbf{V}\cdot N-\int_\Omega \mathcal{C}_i(v_i)\mathbf{Q}\partial_3 \varphi  
		+\int_\Omega \mathbf{F}_{ik} \mathcal{R}^2_{ik} \partial_3 \varphi+\int_\Omega \mathcal{R}_i^1\mathbf{V}_i \partial_3 \varphi \label{spatial}.
	\end{align}\\
	
	\noindent\textbf{Control of the error terms}\\
	Let $\partial_3\varphi\geq c_0>0$ and $D^\alpha = \partial_1^{\alpha_1}\partial_2^{\alpha _2}$. It directly follows from the definition of $\mathcal{C}_i(f)$ and $\mathcal{D}(f)$ that
	\begin{align}
		||\mathcal{C}_i(f)||_0&\leq P(c_0^{-1},|\psi|_4)||f||_4 \label{Cest},\\
		||\mathcal{D}(f)||_0&\leq P(c_0^{-1},||v||_4,|\psi|_4,|\partial_t \psi|_3)(||f||_4+||\partial_t f||_2) \label{Dest}.
	\end{align}
	Then
	\begin{align}
	-\int_\Omega\mathcal{C}_i(v_i)\mathbf{Q}\partial_3 \varphi \leq||\mathcal{C}_i(v_i)||_0 ||\mathbf{Q}||_0||\partial_3 \varphi||_\infty,	
	\end{align}
	where 
	\begin{align}
	||\mathbf{Q}||_0	\leq ||D^\alpha q||_0+||D^\alpha \varphi||_0 ||\partial_3 ^\varphi q||_\infty.
	\end{align}
	Next,
	\begin{align}
	\int_\Omega \mathbf{F}_{ik}\mathcal{R}_{ik}^2\partial_3 \varphi \leq ||\partial_3 \varphi||_\infty ||\mathbf{F}||_0 ||\mathcal{R}^2_{ik}||_0,
	\end{align}
	where 
	\begin{align}
	||\mathcal{R}_{ik}^2||_0	&=||F_{kj}\mathcal{C}_k(v_i)+[D^\alpha,F_{kj}]\partial_k v_i-\mathcal{D}(F_{ij})||_0\\
	&\lesssim ||F_{kj}||_\infty||\mathcal{C}_k(v_i)||_0 +||F_{kj}||_4||v||_4+||\mathcal{D}(F_{ij})||_0.
	\end{align}
	Similarly, 
	\begin{align}	
	\int_\Omega \mathcal{R}_i^1\mathbf{V}_i \partial_3 \varphi\leq||\mathcal{R}_i^1 ||_0||\mathbf{V}_i||_0||\partial_3\varphi||_\infty,
	\end{align}
	where
	\begin{align}
	||\mathcal{R}_i^1||_0&=	||F_{lk}\mathcal{C}_l(F_{ik})+[D^\alpha, F_{lk}] \partial_l^\varphi F_{ik}-\mathcal{D}(v_i)-\mathcal{C}_i(q)||_0\\
	&\lesssim||F_{ik}||_\infty ||\mathcal{C}_l(F_{ik})||_0+||F_{lk}||_4^2+||\mathcal{D}(v_i)||_0+||\mathcal{C}_i(q)||_0.
	\end{align}\\
	
	\noindent\textbf{Control of $-\int_\Sigma \mathbf{Q}(\mathbf{V}\cdot N)$}\\
	Plugging in the higher order kinematic boundary condition (\ref{kbc}), we have 
	\begin{align}
	-\int_\Sigma \mathbf{Q} (\mathbf{V}\cdot N) = -\int_\Sigma \mathbf{Q}(\partial_t^\varphi D^\alpha \psi +\overline{v}\cdot D^\alpha \overline{\partial}\psi -\mathcal{S}_1).
	\end{align}
	For the first term, invoking the boundary condition for pressure (\ref{bcq}),
	\begin{align}
	I:&=-\int_\Sigma \mathbf{Q}\partial_t^\varphi D^\alpha \psi\nonumber \\
	&=\int_\Sigma\sigma D^\alpha\big(\overline{\nabla}\cdot \frac{\overline{\nabla}\psi}{|N|}\big)\partial_t^\varphi D^\alpha \psi
		+\partial_3 qD^\alpha \psi \partial_t^\varphi D^\alpha \psi \nonumber\\
	&=:ST+RT.
	\end{align}
	We first deal with the term contributed by the surface tension. Let $|\beta|=1$, integrating $\overline{\nabla}$ by parts, 
	\begin{align}
	ST&=-\int_\Sigma\sigma D^\alpha\big(\frac{\overline{\nabla}\psi}{|N|}\big)\cdot \partial_t^\varphi \overline{\nabla}D^\alpha \psi \nonumber\\
	&=-\int_\Sigma \sigma \big(\frac{D^\alpha\overline{\nabla}\psi}{|N|}- \frac{\overline{\nabla}\psi\cdot D^\alpha \overline{\nabla}\psi}{|N|^3}\overline{\nabla}\psi \big)\cdot \partial_t^\varphi \overline{\nabla}D^\alpha \psi \nonumber\\
	 &\quad -\int_\Sigma \sigma\big([D^\alpha,\overline{\nabla}\psi,\frac{1}{|N|}] -\overline{\nabla}\psi [D^{\alpha-\beta},\frac{\overline{\nabla}\psi}{|N|^3}]\cdot D^\beta \overline{\nabla}\psi \big)\cdot \partial_t^\varphi \overline{\nabla}D^\alpha \psi \nonumber\\
	 &=:ST_1+ST^R_2 \label{ST}.
	\end{align}
	Observing the symmetry in the first term, we expect it to contribute to energy terms,
	\begin{align}
		ST_1&=-\frac{\sigma}{2}\frac{d}{dt}\int_\Sigma \big( \frac{|D^\alpha \overline{\nabla }\psi|^2}{|N|}-\frac{|\overline{\nabla}\psi\cdot D^\alpha \overline{\nabla}\psi|^2}{|N|^3} \big) \nonumber\\
		&\quad+\frac{\sigma}{2}\int_\Sigma \partial_t(\frac{1}{|N|})|D^\alpha\overline{\nabla}\psi|^2-\partial_t(\frac{1}{|N|^3})|\overline{\nabla}\psi\cdot D^\alpha \overline{\nabla}\psi|^2 \nonumber \\
		&=:ST_{11}+ST_{12}^R \label{ST1}.
	\end{align}

	We do the following calculation to check that $ST_{11}$ indeed give rise to the energy of surface tension, we do the following calculation.
	\begin{align}
		\frac{|D^\alpha \overline{\nabla }\psi|^2}{|N|}-\frac{|\overline{\nabla}\psi\cdot D^\alpha \overline{\nabla}\psi|^2}{|N|^3}
		\geq\frac{|D^\alpha \overline{\nabla}\psi|^2(1+|\overline{\nabla}\psi|^2)-|\overline{\nabla}\psi|^2|D^\alpha\overline{\nabla}\psi|^2}{|N|^3}=\frac{|D^\alpha \overline{\nabla}\psi|^2}{|N|^3}.
	\end{align}
	Since $|N|$ is bounded in short time, it follows that 
	\begin{align}
		\frac{\sigma}{2} \int _\Sigma |D^\alpha \overline{\nabla}\psi|^2 
		\leq \frac{\sigma}{2} \int _\Sigma \big(\frac{|D^\alpha \overline{\nabla }\psi|^2}{|N|}-\frac{|\overline{\nabla}\psi\cdot D^\alpha \overline{\nabla}\psi|^2}{|N|^3}\big) \label{coercive}.
	\end{align}

	Now, it remains to control $ST_{2}^R$ and $ST_{12}^R$. For $ST_{12}^R$, observe that $\partial_t(\frac{1}{|N|}),\partial_t(\frac{1}{|N|^3})$ contributes to $|\partial_t\overline{\nabla}\psi|_\infty$. Although $|D^\alpha \overline{\nabla}\psi|_0^2$ contributes to the top order of $\psi$, it is attached with $\sigma$. Hence, 
	\begin{align}
	ST_{12}^R\leq P(|\overline{\nabla}\psi|_\infty)|\overline{\nabla}\partial_t \psi |_\infty |\sqrt{\sigma}D^\alpha \overline{\nabla}\psi|_0^2 \label{ST12}.
	\end{align}
	For $ST_2^R$, we first deal with 
	\begin{align}ST_{21}^R:=-\int_\Sigma \sigma [D^\alpha,\overline{\nabla}\psi,\frac{1}{|N|}]\cdot \partial_t^\varphi \overline{\nabla}D^\alpha \psi. 	
	\end{align}
	Let $|\gamma|=|\zeta|=1$. Notice that in the commutator $[D^\alpha,\overline{\nabla}\psi,\frac{1}{|N|}]$, the top order terms appears when $D^{\alpha-\gamma}$ lands on $\overline{\nabla}\psi$ or $\frac{1}{|N|}$, more precisely, $ST_{21}^R$ contributes to 
	\begin{align}
	ST_{211}^R&:=-\sigma \int_\Sigma D^{\alpha-\gamma}\overline{\nabla} \psi D^\gamma(\frac{1}{|N|})\cdot \partial_t^\varphi \overline{\nabla}D^\alpha \psi ,\\
	ST_{212}^R&:=-\sigma \int_\Sigma D^{\alpha-\gamma-\zeta}\overline{\nabla} \psi D^{\gamma+\zeta}(\frac{1}{|N|})\cdot \partial_t^\varphi \overline{\nabla}D^\alpha \psi ,	\\
	ST_{213}^R&:=-\sigma \int_\Sigma D^{\gamma}\overline{\nabla} \psi D^{\alpha-\gamma}(\frac{1}{|N|})\cdot \partial_t^\varphi \overline{\nabla}D^\alpha \psi.
	\end{align}
	The first term and the last term have similar structures, the second term is even easier as it only contains the lower order term of the commutator, so we only deal with $R_{211}^R$. Integrating $\overline{\nabla}$ by parts, 
	\begin{align}
		ST_{211}^R&=\sigma \int_\Sigma \overline{\nabla}\cdot(D^{\alpha-\gamma}\overline{\nabla} \psi D^\gamma(\frac{1}{|N|})\big)\partial_t^\varphi D^\alpha\psi \nonumber \\
		&=\sigma \int_\Sigma (D^{\alpha-\gamma}\overline{\nabla}\cdot\overline{\nabla}\psi) D^\gamma(\frac{1}{|N|})\partial_t^\varphi D^\alpha \psi +\sigma \int_\Sigma \big(D^{\alpha -\gamma}\overline{\nabla}\psi \cdot \overline{\nabla}D^\gamma (\frac{1}{|N|})\big)\partial_t^\varphi D^\alpha \psi \nonumber\\
		&\lesssim |\psi|_4|\sqrt{\sigma}\overline{\nabla}D^\alpha\psi|_0|\sqrt{\sigma}\partial_t^\varphi D^\alpha \psi|_0 \label{ST211}.
	\end{align}
	The control of $ ST_{22}^R:=\int_\Sigma \sigma \big([D^{\alpha-\beta},\frac{\overline{\nabla}\psi}{|N|^3}]\cdot D^\beta \overline{\nabla}\psi \big)\big(\overline{\nabla}\psi\cdot \partial_t^\varphi \overline{\nabla}D^\alpha \psi\big)$ is similar by noticing that the highest order term of the commutator occurs when $D^{\alpha-\beta}$ lands on $\frac{\overline{\nabla}\psi}{|N|^3}$ or when $D^{\alpha-\beta-\gamma}$ lands on $D^\beta\overline{\nabla}\psi$,
	\begin{align}
		ST_{22}^R\lesssim P(|\psi|_4)|\sqrt{\sigma}\overline{\nabla}D^\alpha \psi|_0|\sqrt{\sigma}\partial_t^\varphi D^\alpha \psi|_0 \label{ST22}.
	\end{align}
	By (\ref{ST}), (\ref{ST1}), (\ref{ST12}), (\ref{ST211}), (\ref{ST22}), we have
	\begin{align}
	ST+	\frac{\sigma}{2}\frac{d}{dt}\int_\Sigma \big( \frac{|D^\alpha \overline{\nabla }\psi|^2}{|N|}-\frac{|\overline{\nabla}\psi\cdot D^\alpha \overline{\nabla}\psi|^2}{|N|^3} \big) \lesssim P(E(t)).
	\end{align}
	It follows from (\ref{coercive}) that 
	\begin{align}
	\int_0^TST +\frac{\sigma}{2}\int_\Sigma |D^\alpha \overline{\nabla}\psi|^2\lesssim P(E(0))+\int_0^T P(E(t)) \label{STe}.
	\end{align}

	Next, we deal with the term $RT=\int_\Sigma \partial_3 qD^\alpha \psi \partial_t^\varphi D^\alpha \psi $. If we assume the Rayleigh-Taylor sign condition $-\partial_3q\geq c_0>0$ on $\Sigma$, then the term $\int _\Sigma |D^\alpha \psi |^2$ enters the energy. More precisely, by observing the symmetry, we have 
	\begin{align}
		RT&=-\int_\Sigma (-\partial_3 q)D^\alpha \psi \partial_t^\varphi D^\alpha \psi \nonumber\\
		&=-\frac{d}{dt}\int_\Sigma(-\partial_3 q)|D^\alpha \psi|^2 + \int_\Sigma \partial_3 \partial_t q|D^\alpha \psi |^2 \label{RTs1},
	\end{align}
	where 
	\begin{align}
	\int_\Sigma \partial_3 \partial_t q|D^\alpha \psi |^2\leq  |\partial_3\partial_tq|_\infty |D^\alpha \psi |^2_0\label{RTs2}.
	\end{align}
	However, if we drop the Rayleigh-Taylor sign condition, we can only control the term depending on $\sigma$, 
	\begin{align}
	RT&=\int_\Sigma (\partial_3 q)D^\alpha \psi \partial_t^\varphi D^\alpha \psi \leq |\partial_3 q|_\infty	|D^\alpha \psi|_0|\partial_t D^\alpha \psi|_0,\label{RTws1}
	\end{align}
	where 
	\begin{align}
		|\partial_t D^\alpha \psi|_0\leq \sigma^{-\frac{1}{2}}\sqrt{E(t)}.\label{RTws2}
	\end{align}
	
	For the second term, 
	\begin{align}
	II:&=-\int_\Sigma \mathbf{Q}(\overline{v}\cdot D^\alpha \overline{\partial}\psi)\nonumber\\
	&=\int_\Sigma\sigma D^\alpha\big(\overline{\nabla}\cdot \frac{\overline{\nabla}\psi}{|N|}\big)(\overline{v}\cdot D^\alpha \overline{\partial}\psi)+\int_\Sigma\partial_3 qD^\alpha \psi(\overline{v}\cdot D^\alpha \overline{\partial}\psi)\nonumber\\
	&=:II_1+II_2. \label{II}
	\end{align}
	For $II_2$, it is not hard to see the symmetry. Integrating $\overline{\partial}$ by parts,
	\begin{align}
	II_2&=-\int_\Sigma 	\partial_3 q(\overline{v}\cdot D^\alpha \overline{\partial}\psi )D^\alpha \overline{\partial}\psi-\int_\Sigma\overline{\partial}\cdot(\partial_3 q \overline{v})|D^\alpha \psi |^2\nonumber\\
	&=-II_2-\int_\Sigma\overline{\partial}\cdot(\partial_3 q \overline{v})|D^\alpha \psi |^2.
	\end{align}
	It follows that 
	\begin{align}
		II_2=-\frac{1}{2}\int_\Sigma\overline{\partial}\cdot(\partial_3 q \overline{v})|D^\alpha \psi |^2\leq P(||p||_4,||v||_3)|D^\alpha \psi |_0^2 \label{II2}.
	\end{align}
	For $II_1$, there is a symmetric structure of $D^\alpha\overline{\partial}\psi$. Let $|\beta|=1$, we expand $D^\alpha(\frac{\overline{\nabla}\psi}{|N|})$
		\begin{align}
		II_1&=\sigma\int_\Sigma \overline{\nabla}\cdot \big(\frac{D^\alpha \overline{\nabla}\psi}{|N|}		-\frac{\overline{\nabla}\psi\cdot D^\alpha \overline{\nabla}\psi}{|N|^3}\overline{\nabla}\psi\big) (\overline{v}\cdot \overline{\nabla})D^\alpha \psi\nonumber\\
		&\quad +\sigma\int_\Sigma \overline{\nabla}\cdot \big([D^\alpha,\overline{\nabla}\psi,\frac{1}{|N|}] -\overline{\nabla}\psi [D^{\alpha-\beta},\frac{\overline{\nabla}\psi}{|N|^3}]\cdot D^\beta \overline{\nabla}\psi\big)(\overline{v}\cdot \overline{\nabla})D^\alpha \psi\nonumber\\
		&=:II_{11}+II_{12}^R.  \label{II1}
	\end{align}
	The commutator term $II_{12}^R$ can be directly controlled similar to $ST_2^R$, we have 
	\begin{align}
	II_{12}^R\lesssim P(|\psi|_4)|v|_0|\sqrt{\sigma}D^\alpha \overline{\nabla}\psi|^2_0\label{II12}.
	\end{align}
	For $II_{11}$, we first integrate $\overline{\nabla}$ by parts, 
	\begin{align}
	II_{11}&=-\sigma \int_\Sigma \big(\frac{D^\alpha \overline{\nabla}\psi}{|N|}		-\frac{\overline{\nabla}\psi\cdot D^\alpha \overline{\nabla}\psi}{|N|^3}\overline{\nabla}\psi\big) \cdot (\overline{v}\cdot \overline{\nabla})D^\alpha \overline{\nabla}\psi \nonumber \\
	&\quad-\sigma \int_\Sigma \bigg(\big(\frac{D^\alpha \overline{\nabla}\psi}{|N|}		-\frac{\overline{\nabla}\psi\cdot D^\alpha \overline{\nabla}\psi}{|N|^3}\overline{\nabla}\psi\big) \cdot \overline{\nabla} \bigg)\overline{v}\cdot \overline{\nabla}D^\alpha \psi\nonumber\\
	&=:II_{111}+II_{112}^R \label{II11},
	\end{align}
	where $II_{112}^R$ can be directly controlled,
	\begin{align}
	II_{112}^R=	-\sigma \int_\Sigma \bigg(\big(\frac{D^\alpha \overline{\nabla}\psi}{|N|}		-\frac{\overline{\nabla}\psi\cdot D^\alpha \overline{\nabla}\psi}{|N|^3}\overline{\nabla}\psi\big) \cdot \overline{\nabla} \bigg)\overline{v}\cdot \overline{\nabla}D^\alpha \psi\lesssim P(|\psi|_3)|\overline{\nabla} \overline{v}|_\infty|\sqrt{\sigma}D^\alpha \overline{\nabla}\psi|_0^2  \label{II112},
	\end{align}
	then we integrate $(\overline{v}\cdot \overline{\nabla})$ in $II_{111}$ by parts to get symmetric terms,
\begin{align}
II_{111}
&=\sigma \int_\Sigma \big(\frac{(\overline{v}\cdot \overline{\nabla})D^\alpha \overline{\nabla}\psi}{|N|}		-\frac{\overline{\nabla}\psi\cdot (\overline{v}\cdot \overline{\nabla})D^\alpha \overline{\nabla}\psi}{|N|^3}\overline{\nabla}\psi\big) \cdot D^\alpha \overline{\nabla}\psi \nonumber\\
&\quad +\sigma \int_\Sigma \big(\frac{D^\alpha \overline{\nabla}\psi}{|N|}		-\frac{\overline{\nabla}\psi\cdot D^\alpha \overline{\nabla}\psi}{|N|^3}\overline{\nabla}\psi\big) \cdot (\overline{\nabla}\cdot \overline{v})D^\alpha \overline{\nabla}\psi
+ \sigma \int _\Sigma |D^\alpha \overline{\nabla}\psi|^2(\overline{v}\cdot \overline{\nabla})(\frac{1}{|N|})		 \nonumber\\
&\quad -\int_\Sigma\bigg(\frac{(\overline{v}\cdot \overline{\nabla})\overline{\nabla}\psi\cdot D^\alpha \overline{\nabla}\psi}{|N|^3}\overline{\nabla}\psi
-\frac{\overline{\nabla}\psi\cdot D^\alpha \overline{\nabla}\psi}{|N|^3}(\overline{v}\cdot \overline{\nabla})\overline{\nabla}\psi
-(\overline{\nabla}\psi\cdot D^\alpha \overline{\nabla}\psi)(\overline{\nabla}\psi\big)(\overline{v}\cdot \overline{\nabla})(\frac{1}{|N|^3})\bigg) \cdot D^\alpha \overline{\nabla}\psi \nonumber\\
&=-II_{111} +II_{1112}^R+II_{1113}^R+II_{1114}^R.
\end{align}
Hence, 
\begin{align}
II_{111}=\frac{1}{2}(II_{1112}^R+II_{1113}^R+II_{1114}^R) \lesssim P(|\psi|_4,||v||_3)|\sqrt{\sigma}D^\alpha \overline{\nabla}\psi|_0^2\label{II111},
\end{align}
which closes the control of $II$. By (\ref{II}), (\ref{II1}), (\ref{II2}), (\ref{II1}), (\ref{II12}), (\ref{II11}), (\ref{II112}), (\ref{II111}), we have 
\begin{align}
\int_0^T II\lesssim P(E(0))+\int_0^TP(E(t))	\label{IIes}.
\end{align}

It remains to control the last term 
\begin{align}
III:&=	\int_\Sigma \mathbf{Q}\mathcal{S}_1\nonumber\\
&=\int_\Sigma D^\alpha q\mathcal{S}_1- \int_\Sigma D^\alpha \psi \partial_3 q\mathcal{S}_1\nonumber\\
&=:III_1+III_2\label{III}.
\end{align}
Since $\mathcal{S}_1$ consists of only lower order terms, we have 
\begin{align}
 III_2\leq |\partial_3 q|_\infty|D^\alpha \psi|_0|\mathcal{S}_1|_0\lesssim	
 P(|\psi|_4,||v||_4)|\partial_3 q|_\infty|D^\alpha \psi|_0.\label{III2}
 \end{align}
For $III_1$, since $D^\alpha q$ has 4 derivatives on the boundary, we have to remove at least half of the derivative of $q$ to close the estimate. Let $|\beta_1'|=3, |\beta_2'|=1$,
\begin{align}
III_1&=\int_\Sigma D^\alpha qD^{\beta_1'}\overline{v}\cdot D^{\beta_2'}\overline{\partial}\psi+\int_\Sigma D^\alpha q\bigg(D^\alpha \psi \partial_3 v\cdot N-\sum_{\substack {|\beta_1|+|\beta_2|=4 \\ \beta_1>0,\beta_2>0\\ \beta_1\neq 3, \beta_2\neq 1}}D^{\beta_1} \overline{v}\cdot D^{\beta_2}\overline{\partial} \psi \bigg)\nonumber\\
&=III_{11}+III_{12}\label{III1}.
\end{align}
For $III_{12}$, we plug in the boundary condition of $q$ , then integrate $\overline{\nabla}$ by parts to get 
\begin{align}
III_{12}&=-\sigma \int _\Sigma D^\alpha (\overline{\nabla}\cdot \frac{\overline{\nabla}\psi}{|N|})\bigg(D^\alpha \psi \partial_3 v\cdot N-\sum_{\substack {|\beta_1|+|\beta_2|=4 \\ \beta_1>0,\beta_2>0\\ \beta_1\neq 3, \beta_2\neq 1}}D^{\beta_1} \overline{v}\cdot D^{\beta_2}\overline{\partial} \psi \bigg)\nonumber\\
&=\sigma \int_\Sigma D^\alpha(\frac{\overline{\nabla}\psi}{|N|})\cdot \overline{\nabla} \bigg(D^\alpha \psi \partial_3 v\cdot N-\sum_{\substack {|\beta_1|+|\beta_2|=4 \\ \beta_1>0,\beta_2>0\\ \beta_1\neq 3, \beta_2\neq 1}}D^{\beta_1} \overline{v}\cdot D^{\beta_2}\overline{\partial} \psi \bigg)\nonumber\\
&\lesssim P(|\psi|_4,||v||_4,|\sqrt{\sigma}D^\alpha \overline {\nabla} \psi |_0)\label{III12}.
\end{align}
$III_{11}$ can be directly controlled by
\begin{align}
III_{11}\lesssim |\psi|_4|D^\alpha q|_{-\frac{1}{2}}|D^{\beta_1'}v|_{\frac{1}{2}}\lesssim |\psi|_4 ||q||_4 ||v||_4\label{III11}.
\end{align}

Combining (\ref{spatial}), (\ref{STe}), (\ref{RTs1}), (\ref{RTs2}), (\ref{RTws1}), (\ref{RTws2}), (\ref{IIes}), (\ref{III}), (\ref{III2}), (\ref{III1}),(\ref{III11}), (\ref{III12}), we have 
\begin{align}
	||\mathbf{V}_i||_0^2 +||\mathbf{F}_{ik}||_0^2+|\sqrt{\sigma}D^\alpha \overline{\nabla}\psi|_0^2\lesssim P(E(0))+\int_0^TP(E(t)).
\end{align}
Furthermore, if the Rayleigh-Taylor sign condition is assumed, the term $|D^\alpha \psi|_0^2$ enters the energy and we have the $\sigma-$uniform estimate
\begin{align}
	||\mathbf{V}_i||_0^2 +||\mathbf{F}_{ik}||_0^2+|\sqrt{\sigma}D^\alpha \overline{\nabla}\psi|_0^2+|D^\alpha \psi|_0^2\lesssim P(E(0))+\int_0^TP(E(t)).
\end{align}
By (\ref{AGU}), we can replace $\mathbf{V}$ and $\mathbf{F}$ by $D^\alpha v$ and $D^\alpha F$ respectively. Hence, we conclude the proof of \textbf{Theorem \ref{spatial-estimate}}.

The case when there is at least one spatial derivative in $D^\alpha$, i.e. $0<\alpha_0\leq 3$ can be studied with the same analysis as above.
\subsubsection{Energy estimate with time derivatives}
 In this section, we focus on the fully times differentiated equations, i.e $D^\alpha=\partial_t^4$.

\begin{theorem}\label{time-estimate}
	There exists $T>0$ such that 
	\begin{align}
	||\partial_t^4 v_i||_0^2	+||\partial_t^4  F_{ik}||_0^2+|\sqrt{\sigma}\partial_t^4 \overline{\nabla}\psi|_0^2\lesssim P(E(0))+\int_0^T P(E(t)).
	\end{align}
	Furthermore, if the Rayleigh-Taylor sign condition $\-\partial_3 q\geq c_0 >0$ is assumed on $\Sigma$, the term $|\partial_t^4\psi|_0^2$ enters the energy and the estimate become $\sigma$-uniform,
	\begin{align}
		||\partial_t^4 v_i||_0^2	+||\partial_t^4  F_{ik}||_0^2+|\sqrt{\sigma}\partial_t^4 \overline{\nabla}\psi|_0^2+|\partial_t^4\psi|_0^2\lesssim P(E(0))+\int_0^T P(E(t)).
	\end{align}

	\end{theorem}
Parallel to the case with spatial derivatives (\ref{spatial}), we have 
\begin{align}
		\frac{1}{2}\frac{d}{dt}(\int_\Omega |\mathbf{V}_i|^2\partial_3 \varphi +\int_\Omega |\mathbf{F}_{ik}|^2 \partial_3 \varphi)=
			-\int_\Sigma \mathbf{Q}\mathbf{V}\cdot N-\int_\Omega \mathcal{C}_i(v_i)\mathbf{Q}\partial_3 \varphi  
		+\int_\Omega \mathbf{F}_{ik} \mathcal{R}^2_{ik} \partial_3 \varphi+\int_\Omega \mathcal{R}_i^1\mathbf{V}_i \partial_3 \varphi. \label{time}
	\end{align}
The Alinhac good unknowns reads
\begin{align}
	\mathbf{V}=\partial_t^4 v - \partial_t^4 \varphi\partial_3^\varphi v,\quad
	\mathbf{F}=\partial_t^4 F - \partial_t^4 \varphi\partial_3^\varphi F,\quad
	\mathbf{Q}=\partial_t^4 q - \partial_t^4 \varphi\partial_3^\varphi q
	\end{align}
	with the following properties, for any function $g$ and $i=1,2,3,$
\begin{align}
\partial_t^4 \partial_i^\varphi g&=\partial_i^\varphi \mathbf{G}+\mathcal{C}_i(g),\\
\partial_t^4 D_t^\varphi g&=D_t^\varphi \mathbf{G}+\mathcal{D}(g),
\end{align}
where, for $\tau=1,2,$ 
	\begin{align}
		\mathcal{C}_\tau(g)&=\partial_t^4\varphi\partial_\tau^\varphi\partial_3^\varphi g-[\partial_t^4,\frac{\partial_\tau \varphi}{\partial_3 \varphi },\partial_3 g]-\partial_3 g[\partial_t^4,\partial_\tau \varphi,\frac{1}{\partial_3 \varphi}]+\partial_3 g \partial _\tau \varphi[\partial_t^3,\frac{1}{(\partial_3 \varphi)^2}]\partial_t\partial_3 \varphi,\\
		\mathcal{C}_3(g)&=\partial_t^4 \varphi (\partial_3 ^\varphi)^2 g+[\partial_t^4,\frac{1}{\partial_3 \varphi},\partial_3 g]- \partial_3 g [\partial_t^3,\frac{1}{(\partial_3 \varphi)^2}]\partial_t \partial_3 \varphi,\\
		\mathcal{D}(g)&=\partial_t^4 \varphi D_t^\varphi \partial_3^\varphi g +[\partial_t^4, \overline{v}]\cdot \overline{\nabla}g+[\partial_t^4, \frac{1}{\partial_3 \varphi}(v\cdot \mathbf N-\partial_t \varphi),\partial_3 g]+[\partial_t^4, \frac{1}{\partial_3 \varphi},v\cdot \mathbf{N}-\partial_t \varphi]\partial_3 g\nonumber\\
			&\quad-(v\cdot \mathbf{N} -\partial_t \varphi)\partial_3 g[\partial_t^3,\frac{1}{(\partial_3 \varphi)^2}]\partial_t\partial_3 \varphi + \frac{1}{\partial_3 \varphi}\partial_3 g[\partial_t^4,v]\mathbf{N}.
	\end{align}
Also, the remaining terms read
	\begin{align}
	\mathcal{R}_i^1&=F_{lk}\mathcal{C}_l(F_{ik})+[\partial_t^4, F_{lk}] \partial_l^\varphi F_{ik}-\mathcal{D}(v_i)-\mathcal{C}_i(q)\label{R1t},\\
	\mathcal{R}_{ij}^2&=F_{kj}\mathcal{C}_k(v_i)+[\partial_t^4,F_{kj}]\partial_k v_i-\mathcal{D}(F_{ij})	\label{R2t}.
	\end{align}
	
\noindent\textbf{Control of the error terms}\\
Let $\partial_3 \varphi \geq c_0>0$, then 
\begin{align}
||\mathcal{C}_i(g)||_0&\leq P\bigg(c_0^{-1},|\overline{\nabla}\psi|_\infty, \sum_{k=1}^3 |\overline{\nabla}\partial_t^k\psi|_{3-k}, |\partial_t^4\psi |_0\bigg) \cdot \bigg(||\partial g||_\infty +||\partial^2 g||_\infty + \sum _{k=1}^3 ||\partial_t^k g||_{4-k}\bigg),\\
||\mathcal{D}(g)||_0&\leq P\bigg(c_0^{-1},\sum_{k=0}^3||\partial_t^k v||_0,|\overline{\nabla}\psi|_\infty, \sum_{k=1}^3 |\overline{\nabla}\partial_t^k\psi|_{3-k}, |\partial_t^4\psi |_0\bigg) \cdot \bigg(||\partial g||_\infty +||\partial^2 g||_\infty + \sum _{k=1}^3 ||\partial_t^k g||_{4-k}\bigg),
\end{align}
where $|\partial_t^4 \psi|_0$ can be controlled by taking $\partial_t^3$ on the kinematic boundary condition. More precisely,
\begin{align}
\partial_t^4\psi
&=\partial_t ^3 (v\cdot N)\nonumber\\
&=-\overline{v}\cdot \overline{\partial}\partial_t ^3 \psi- \partial_t^3 \overline{v}\cdot \overline{\partial}\psi+\partial_t^3 v_3 -[\partial_t^3,\overline{v},\overline{\partial}\psi ]\nonumber\\
&=-\overline{v}\cdot \overline{\partial}\partial_t ^3 \psi-\partial_t^3v \cdot N -[\partial_t^3,\overline{v},\overline{\partial}\psi ]. \label{pt4psi}
\end{align}
It follows that 
\begin{align}
|\partial_t^4\psi|_0\leq P\big(\sum _{k=0}^3|\overline{\nabla }\partial_t^k \psi|_{3                      -k}, \sum_{k=0}^3||\partial_t^k v||_{3-k} \big).
\end{align}
Next,
\begin{align}
||\mathcal{R}_i^1||_0 &\leq P\bigg(\sum_{k=0}^4||\partial_t ^k F||_{4-k}\bigg)\cdot \bigg(  ||\mathcal{C}_l(F_{ik})||_0 + ||\mathcal{D}(v_i)||_0 +\mathcal{C}_i(q)||_0\bigg),\\
||\mathcal{R}_{ij}^2||_0 &\leq P\bigg(\sum_{k=0}^4||\partial_t ^k F||_{4-k}, \sum_{k=0}^3||\partial_t^k v||_{4-k}\bigg)\cdot \bigg(  ||\mathcal{C}_l(v_i)||_0 + ||\mathcal{D}(F_{ij})||_0\bigg).
\end{align}
Then, we can directly control the error terms
\begin{align}
\int_\Omega \mathbf{F}_{ik} \mathcal{R}^2_{ik} \partial_3 \varphi&\leq ||\partial_3 \varphi||_\infty||\mathbf{F}_{ik}||_0||\mathcal{R}_{ik}^2||_0,\\
\int_\Omega \mathcal{R}_i^1\mathbf{V}_i \partial_3 \varphi 	&\leq ||\partial_3 \varphi ||_\infty ||\mathbf{V}_i||_0||\mathcal{R}_i^1||_0.
\end{align}
As for 
\begin{align}
-\int _\Omega \mathcal{C}_i(v_i)\mathbf{Q}\partial_3 \varphi	= -\int_\Omega \partial_t^4q \mathcal{C}_i(v_i) \partial_3\varphi
+\int_\Omega \partial_t^4 \psi \partial_3 q\mathcal{C}_i(v_i)\partial_3 \varphi,
\end{align}
we can only control the second term
\begin{align}
\int_\Omega \partial_t^4 \psi \partial_3 q\mathcal{C}_i(v_i)\partial_3 \varphi\leq |\partial_t^4 \psi|_0|\partial_3q|_\infty||\mathcal{C}_i(v_i)||_0|\partial_3\varphi|_\infty,
\end{align}
and $-\int_\Omega \partial_t^4q \mathcal{C}_i(v_i) \partial_3\varphi
 $ cannot be directly controlled as above, since we do not have the control of $||\partial_t^4 q||_0$. However, it can be cancelled later in the control of $\int_\Sigma \mathbf{Q}\mathcal{S}_1$.\\

\noindent\textbf{Control of $-\int_\Sigma \mathbf{Q}\mathbf{V}\cdot N$  }\\
Plugging in $D^\alpha =\partial_t^4$ to (\ref{kbc}), (\ref{S1}), we have the time differentiated kinematic boundary condition
\begin{align}
\partial_t^5\psi = \mathbf{V}\cdot N- \overline{v}\cdot \partial_t^4\overline{\nabla}\psi +\mathcal{S}^*_1, \label{kbct}
\end{align}
where 
\begin{align}
\mathcal{S}^*_1=\partial_t^4 \psi \partial_3 v\cdot N-[\partial_t^4, \overline{v},\overline{\partial}\psi].
\end{align}
Also, we have the boundary condition for $\mathbf{Q}$ on $\Sigma$,
\begin{align}
\mathbf{Q}|_{\Sigma}=-\sigma \partial_t^4\big(\overline{\nabla}\cdot \frac{\overline{\nabla}\psi }{|N|}\big)-\partial_t^4\psi \partial_3 q.
\end{align}
Now, we plug in the kinematic boundary condition to get 
\begin{align}
	-\int_\Sigma \mathbf{Q}\big(\mathbf{V}\cdot N\big)&= -\int_\Sigma \mathbf{Q} \partial_t^5\psi -\int_\Sigma \mathbf{Q}\big(\overline{v}\cdot \partial_t^4\overline{\nabla}\psi\big) +\int_\Sigma \mathbf{Q}\mathcal{S}^*_1\nonumber\\
	&=:I^*+II^*+III^*.
\end{align}

For $I^*$, plugging in the boundary condition, 
\begin{align}
I^*&=\int_\Sigma \sigma \partial_t^4\big(\overline{\nabla}\cdot \frac{\overline{\nabla}\psi }{|N|}\big)\partial_t^5\psi+\int_\Sigma  \partial_3 q\partial_t^4\psi\partial_t^5\psi\nonumber\\
&=:ST^*+RT^*.
\end{align}
The control of $ST^*$ is parallel to the estimate of $ST$ replacing $D^\alpha$ by $\partial_t^4$, we have 
\begin{align}
\int_0^T ST^*+\frac{\sigma}{2}\int_\Sigma |\partial_t^4 \overline{\nabla}\psi|_0^2 \lesssim P(E(0))+\int_0^TP(E(t))	 \label{ST*}.
\end{align}
For $RT^*$, if we assume the Rayleigh-Taylor sign condition $-\partial_3 q\geq c_0 >0$, then we can handle $RT^*$ similarly to $RT$ by generating an energy term $|\partial_t^4 \psi|_0^2 $. However, if we drop the sign condition, we will need some new estimate to bound $RT^*$, since it has a top order term $\partial_t^5\psi $ which is not in the energy.

To study $RT^*$, we use the kinematic boundary condition,
\begin{align}
\partial_t^5\psi = \partial_t^4(v\cdot N)=
-(\overline{v}\cdot \overline{\partial})\partial_t ^4 \psi+\partial_t^4 v \cdot N -[\partial_t^4,\overline{v},\overline{\partial}\psi ].
\end{align}
It follows that 
\begin{align}
RT^*&=-\int_\Sigma \partial_3 q\partial_t^4\psi(\overline{v}\cdot \overline{\partial})\partial_t ^4 \psi+
\int_\Sigma \partial_3 q\partial_t^4\psi\partial_t^4 v \cdot N -
\int_\Sigma \partial_3 q\partial_t^4\psi[\partial_t^4,\overline{v},\overline{\partial}\psi ]\\
&=:RT^*_1+RT_2^*+RT^*_3 \label{RT^*}.
\end{align}
Note that we do not need to estimate $RT_1^*$ and $RT^*_3$, since they will be cancelled by parts of $II^*$ and $III^*$. As for $RT^*_2$, the term $\partial_t^4 v$ has the top order of $v$ on the boundary, which we certainly cannot directly control. We first plug in (\ref{pt4psi}) for $\partial_t^4\psi$ and use Green's formula to drag the term from boundary to the domain $\Omega$,
\begin{align}
RT_2^*&=	
\int_\Sigma \partial_3 q\big(-\overline{v}\cdot \overline{\partial}\partial_t ^3 \psi-\partial_t^3v \cdot N -[\partial_t^3,\overline{v},\overline{\partial}\psi ] \big)\partial_t^4 v \cdot N\nonumber\\
&=\int_\Omega \partial_3 \bigg( \partial_3 q\big(-\overline{v}\cdot \overline{\partial}\partial_t ^3 \varphi-\partial_t^3v \cdot \mathbf{N} -[\partial_t^3,\overline{v},\overline{\partial}\varphi ] \big)\partial_t^4 v \cdot \mathbf{N}\bigg)\nonumber\\
&=\int_\Omega \partial_3 q\big(-\overline{v}\cdot \overline{\partial}\partial_t ^3 \varphi-\partial_t^3v \cdot \mathbf{N} -[\partial_t^3,\overline{v},\overline{\partial}\varphi ] \big)\partial_3\partial_t^4 v \cdot \mathbf{N}\nonumber\\
&\quad+\int_\Omega \partial_3 q\partial_3\big(-\overline{v}\cdot \overline{\partial}\partial_t ^3 \varphi-\partial_t^3v \cdot \mathbf{N} -[\partial_t^3,\overline{v},\overline{\partial}\varphi ] \big)\partial_t^4 v \cdot \mathbf{N}\nonumber\\
&\quad+\int_\Omega \partial_3^2 q\big(-\overline{v}\cdot \overline{\partial}\partial_t ^3 \varphi-\partial_t^3v \cdot \mathbf{N} -[\partial_t^3,\overline{v},\overline{\partial}\varphi ] \big)\partial_t^4 v \cdot \mathbf{N}\nonumber\\
&\quad+\int_\Omega \partial_3 q\big(-\overline{v}\cdot \overline{\partial}\partial_t ^3 \varphi-\partial_t^3v \cdot \mathbf{N} -[\partial_t^3,\overline{v},\overline{\partial}\varphi ] \big)\partial_t^4 v \cdot \partial_3\mathbf{N} \nonumber\\
&=RT^*_{21}+RT^*_{22}+RT^*_{23}+RT^*_{24} \label{RT^*2},
\end{align}
where $RT^*_{22}, RT^*_{23}, RT^*_{24}$ can be directly controlled
\begin{align}
RT^*_{22}&\leq P\bigg(\sum_{k=0}^3||\partial_t^kv||_{4-k},\sum_{k=0}^3|\partial_t^k\overline{\nabla}\psi|_{4-k}\bigg)||\partial_3q||_\infty	|\overline{\nabla}\psi|_\infty \label{RT*22},\\
RT^*_{23}&\leq P\bigg(\sum_{k=0}^3||\partial_t^kv||_{3-k},\sum_{k=0}^3|\partial_t^k\overline{\nabla}\psi|_{3-k}\bigg)||\partial^2_3q||_\infty	|\overline{\nabla}\psi|_\infty\label{RT*23} ,\\
RT^*_{24}&\leq P\bigg(\sum_{k=0}^3||\partial_t^kv||_{3-k},\sum_{k=0}^3|\partial_t^k\overline{\nabla}\psi|_{3-k}\bigg)||\partial_3q||_\infty	|\partial_3\overline{\nabla}\psi|_\infty \label{RT*24}.
\end{align}
For $RT^*_{21}$, we invoke the time integral and integrate $\partial_t$ by parts,
\begin{align}
\int_0^TRT^*_{21}&=-\int_0^T	\int_\Omega \partial_t\bigg(\partial_3 q\big(-\overline{v}\cdot \overline{\partial}\partial_t ^3 \varphi-\partial_t^3v \cdot \mathbf{N} -[\partial_t^3,\overline{v},\overline{\partial}\varphi ] \big)\mathbf{N}\bigg)\cdot\partial_3\partial_t^3 v \nonumber\\
&\quad+\bigg(\int_\Omega \partial_3 q\big(-\overline{v}\cdot \overline{\partial}\partial_t ^3 \varphi-\partial_t^3v \cdot \mathbf{N} -[\partial_t^3,\overline{v},\overline{\partial}\varphi ] \big)\partial_3\partial_t^3 v \cdot \mathbf{N}\bigg)\bigg|_0^T\nonumber\\
&=J_1+J_2,\label{RT^*21}
\end{align}
where
\begin{align}
J_1&\leq \int_0^T \bigg|\bigg|\partial_t\bigg(\partial_3 q\big(-\overline{v}\cdot \overline{\partial}\partial_t ^3 \varphi-\partial_t^3v \cdot \mathbf{N} -[\partial_t^3,\overline{v},\overline{\partial}\varphi ] \big)\mathbf{N}\bigg)\bigg|\bigg|_0||\partial_3\partial_t^3 v||_0	\nonumber\\
&\leq  \int_0^T P\bigg(\sum_{k=0}^4||\partial_t^kv||_{4-k},\sum_{k=0}^4|\partial_t^k\overline{\nabla}\psi|_{4-k},||\partial_3 q||_\infty, ||\partial_t\partial_3 q||_\infty \bigg) \label{J1}.
\end{align}
As for $J_2$, the only trouble term is $||\partial_3 \partial_t^3 v||_0 $, which we can make absorbed to the left-hand side by Cauchy inequality, 
\begin{align}
J_2\leq P(E(0))+ \epsilon ||\partial_3\partial_t^3 v||_0^2 +C_\epsilon P\bigg(\sum_{k=0}^3||\partial_t^kv||_{3-k},\sum_{k=0}^3|\partial_t^k\overline{\nabla}\psi|_{3-k}\bigg)||\partial_3q||_\infty	|\overline{\nabla}\psi|_\infty.\label{J2}
\end{align}
Combining (\ref{ST*}), (\ref{RT*22}),  (\ref{RT*23}), (\ref{RT*24}), (\ref{J1}), (\ref{J2}), we get the estimate of $I^*$,
\begin{align}
\int_0^T I^*+\frac{\sigma}{2}|\overline{\nabla}\partial_t^4\psi|_0^2 \lesssim \epsilon ||\partial_3\partial_t^3 v||_0^2
+P(E(0))+ \int_0^T P(E(t)) \label{I*est}.
\end{align}

Next, we expand $II^*+III^*$,
\begin{align}
II^*+III^*=	&-\int _\Sigma \partial_t^4q(\overline{v}\cdot \partial_t^4\overline{\nabla}\psi)
+\int _\Sigma \partial_t^4\psi \partial_3 q(\overline{v}\cdot \partial_t^4\overline{\nabla}\psi)\nonumber\\
&+\int_\Sigma \partial_t^4 q \mathcal{S}_1^* 
- \int_\Sigma |\partial_t^4 \psi|^2 \partial_3 q  (\partial_3 v\cdot N)
+\int_\Sigma \partial_t^4 \psi \partial_3 q[\partial_t^4,\overline{v},\overline{\partial}\psi] \label{II*+III*},
\end{align}
where the second term cancels $RT^*_1$ and the fifth term cancels $RT^*_3$. The first term can be controlled exactly as the analysis of $II_{11}$ replacing $D^\alpha$ by $\partial_t^4$ and the fourth term can be directly controlled as follows,
\begin{align}
	- \int_\Sigma |\partial_t^4 \psi|^2 \partial_3 q  (\partial_3 v\cdot N)\leq |\partial_t^4\psi |_0^2|\partial _3 q |_\infty |\partial_3 v|_\infty |N|_\infty.\label{Final1}
\end{align}
It remains to control the third term,
\begin{align}
\int_\Sigma \partial_t^4 q\mathcal{S}_1^*= 	\int_\Sigma \partial_t^4 q \partial_t^4 \psi \partial_3 v\cdot N-4\int_\Sigma \partial_t^4 q (\partial_t^3 \overline{v}\cdot \overline{\partial}\partial_t\psi)-\int _\Sigma \partial_t^4q \sum _{k=1}^2 \begin{pmatrix} 4 \\ k \end{pmatrix} \partial_t^kv\cdot \partial_t^{4-k}\overline{\partial}\psi.
\end{align}
For the first term, invoking the boundary condition and integration by parts,
\begin{align}
	\int_\Sigma \partial_t^4 q \partial_t ^4\psi \partial_3 v \cdot N 
	&=-\int _\Sigma \sigma \partial_t^4 \big(\overline{\nabla}\cdot \frac{\overline{\nabla}\psi}{|N|}\big)\partial_t ^4 \psi \partial_3 v\cdot N \nonumber\\
	&=\int _\Sigma \sigma \partial_t^4 \big( \frac{\overline{\nabla}\psi}{|N|}\big)\cdot \overline{\nabla}\big(\partial_t ^4 \psi \partial_3 v\cdot N\big)\nonumber\\
	&\leq P\big(|\sqrt{\sigma}\partial_t^4 \overline{\nabla}\psi|_0, |\partial_t ^4 \psi|_0,|\partial_3 \overline{\nabla}v|_\infty,|\partial_3 v|_\infty,|\overline{\nabla}^2\psi|_\infty,|\overline{\nabla}\psi|_\infty \big)\label{Final2}.
\end{align}
Similarly, for the last term,
\begin{align}
	&-\int _\Sigma \partial_t^4q \sum _{k=1}^2 \begin{pmatrix} 4 \\ k \end{pmatrix} \partial_t^kv\cdot \partial_t^{4-k}\overline{\partial}\psi\nonumber\\
	&=
	\int _\Sigma \sigma \partial_t^4 \big(\overline{\nabla}\cdot \frac{\overline{\nabla}\psi}{|N|}\big)\sum _{k=1}^2 \begin{pmatrix} 4 \\ k \end{pmatrix} \partial_t^kv\cdot \partial_t^{4-k}\overline{\partial}\psi\nonumber\\
	&=-\int _\Sigma \sigma \partial_t^4 \big( \frac{\overline{\nabla}\psi}{|N|}\big)\cdot\overline{\nabla}\sum _{k=1}^2 \begin{pmatrix} 4 \\ k \end{pmatrix} \partial_t^kv\cdot \partial_t^{4-k}\overline{\partial}\psi\nonumber\\
	&\leq P\big(|\partial_t^2 \overline{\nabla}v|_0,|\partial_t^2 v|_\infty,|\partial_t^2 \overline{\partial}^2\psi|_0,|\partial_t^2 \overline{\partial}\psi|_\infty ,|\partial_t\overline{\nabla} v|_\infty,|\partial_t v|_\infty,|\sqrt{\sigma}\overline{\nabla}^2\partial_t^3\psi|_0,|\overline{\nabla}\partial_t^3\psi|_0 \big)|\sqrt{\sigma}\partial_t^4 \overline{\nabla}\psi|_0\label{Final3}.
\end{align}
As for the second term
\begin{align}
J^*_1:=-4\int_\Sigma \partial_t^4 q (\partial_t^3 \overline{v}\cdot \overline{\partial}\partial_t\psi),
\end{align}
we seek for cancellation from the interior term $-\int_\Omega \partial_t^4q \mathcal{C}_i(v_i) \partial_3\varphi$. Note that the possible trouble terms in $-\int_\Omega \partial_t^4q \mathcal{C}_i(v_i) \partial_3\varphi$ are 
\begin{align}
J_2^*&=-\int_\Omega \partial_t^4q \partial_t^4 \varphi \partial_i^\varphi \partial_3^\varphi v_i \partial_3\varphi, \qquad i=1,2,3,\\
J_3^*&=4\int_\Omega \partial_t^4q \partial_t\big(\frac{\partial_\tau\varphi}{\partial_3 \varphi})\partial_t^3\partial_3 v_\tau \partial_3\varphi, \qquad \tau=1,2,\\
J_4^*&=-4\int_\Omega \partial_t^4q \partial_t\big(\frac{1}{\partial_3\varphi}\big)\partial_t^3 \partial_3 v_3 \partial_3\varphi,	
\end{align}
while the remaining terms are contributed by the lower order terms of $\mathcal{C}_i(v_i)$ and can be controlled by integrating $\partial_t$ by parts under time integral.
Due to the divergence free of $v$, $J_2^*$ vanishes,
\begin{align}
J_2^*=-\int_\Omega \partial_t^4q \partial_t^4 \varphi\partial_3^\varphi (\nabla^\varphi \cdot v) \partial_3\varphi=0.
\end{align}
Next, we expand $J_3^*$,
\begin{align}
J_3^* &=	4\int_\Omega \partial_t^4 q(\partial_t \overline{\partial} \varphi \cdot  \partial_t^3 \partial_3 v) -4 \int_\Omega \partial_t^4 q \frac{\partial_\tau \varphi \partial_t\partial_3 \varphi}{\partial_3 \varphi}\partial_t^3\partial_3 v_\tau \nonumber\\
&=:J_{31}^*+J_{32}^*.
\end{align}
We can easily see the similarity between $J_1^*$ and $J_{31}^*$. After integrating $\partial_3$ by parts in $J_{31}^*$, the boundary term exactly cancels $J_1^*$.
\begin{align}
J_1^*+J_{31}^*=4	\int_\Omega \partial_t^4 \partial_3 q(\partial_t \overline{\partial} \varphi \cdot  \partial_t^3 v) +\partial_t^4  q(\partial_t \partial_3\overline{\partial} \varphi \cdot  \partial_t^3 v),
\end{align}
which can be handled by integrating $\partial_t$ by parts
\begin{align}	
&\int_0^T\int_\Omega \partial_t^4 \partial_3 q(\partial_t \overline{\partial} \varphi \cdot  \partial_t^3 v) +\partial_t^4  q(\partial_t \partial_3\overline{\partial} \varphi \cdot  \partial_t^3 v)\nonumber\\
&=-\int_0^T \int_\Omega \bigg(\partial_t^3 \partial_3 q \partial_t (\partial_t \overline{\partial} \varphi \cdot  \partial_t^3 v) +\partial_t^3 q\partial_t(\partial_t \partial_3\overline{\partial} \varphi \cdot  \partial_t^3 v)\bigg)+\bigg(\int_\Omega \partial_t^3 \partial_3 q  (\partial_t \overline{\partial} \varphi \cdot  \partial_t^3 v+\partial_t^3 q(\partial_t \partial_3\overline{\partial} \varphi \cdot  \partial_t^3 v)) \bigg)\bigg|_0^T \nonumber\\
&\leq  \epsilon ||\partial_t^3 \partial_3 q||_0^2 +P(E(0))+\int_0^TP(E(t)) \label{Final5}.
\end{align}
Finally, it remains to control $J_{32}^*$ and $J_4^*$. Recall that $\partial_\tau^\varphi=\partial_\tau -\partial_\tau \varphi \frac{1}{\partial_3\varphi}\partial_3 $ and $\partial_3 ^\varphi =\frac{1}{\partial_3 \varphi}\partial_3$, so
\begin{align}
J^*_4+J_{32}^*&=4\int _\Omega 	 \partial_t^4 q \partial_t \partial_3 \varphi\partial_3^\varphi \partial_t^3  v_3 
+4 \int_\Omega \partial_t^4 q   \partial_t\partial_3 \varphi\big(\partial_\tau^\varphi \partial_t^3 v_\tau 
- \partial_\tau \partial_t^3 v_\tau\big)\nonumber\\
&=4\int _\Omega 	 \partial_t^4 q \partial_t \partial_3 \varphi\big(\nabla^\varphi\cdot \partial_t ^3  v\big) 
-4\int_\Omega \partial_t^4 q  \varphi \partial_t\partial_3 \varphi\partial_\tau \partial_t^3 v_\tau\nonumber\\
&=4\int _\Omega 	 \partial_t^4 q \partial_t \partial_3 \varphi\partial_t ^3\big(\nabla^\varphi\cdot   v\big)  
-4\int _\Omega 	 \partial_t^4 q \partial_t \partial_3 \varphi\big(\big[\partial_t ^3,\nabla^\varphi\big]\cdot   v\big)
-4\int_\Omega \partial_t^4 q  \partial_t\partial_3 \varphi\partial_\tau \partial_t^3 v_\tau,
\end{align}
where the first term vanishes due to $\nabla^\varphi \cdot v=0$. The second term can be controlled after integration $\partial_t$ by parts since $[\partial_t ^3,\nabla^\varphi\big]\cdot   v$ only contains up to 2 time derivatives of $v$ and 3 time derivatives of $\varphi$. For the last term, after integration $\partial_t$ by parts and integration $\partial_\tau$ by parts, it contributes to the top order term,
\begin{align}
	-4\int_0^T\int_\Omega \partial_t^3 \partial_\tau q  \partial_t\partial_3 \varphi \partial_t^4 v_\tau
	-\bigg(4\int_\Omega \partial_t^3 q   \partial_t\partial_3 \varphi\partial_\tau \partial_t^3 v_\tau\bigg)\bigg|_0^T \leq \epsilon ||\partial_\tau \partial_t^3 v||_0^2+ P(E(0))+ \int_0^T P(E(t)). \label{Final6}
\end{align}
Combining (\ref{I*est}), (\ref{II*+III*}), (\ref{Final1}), (\ref{Final2}),(\ref{Final3}), (\ref{Final5}), (\ref{Final6}), after choosing $\epsilon$ carefully, we conclude the estimate 
\begin{align}
||\mathbf{V}_i||_0^2+||\mathbf{F}_{ik}||_0^2+|\sqrt{\sigma}\partial_t^4 \overline{\nabla}\psi|_0^2 	\leq P(E(0))+\int_0^TP(E(t)).
\end{align}

\section{Zero surface tension limit}
In this section, we study the behaviour of the solution of (\ref{phisys})-(\ref{phibc}) as the surface tension coefficient $\sigma$ tends to 0 and thus show that the solution of (\ref{phisys})-(\ref{phibc}) can be passed to the zero surface tension limit. Consider the following incompressible elastodynamics without surface tension reformulated in graphic coordinates:
	\begin{equation}\label{limitsys}
		\begin{cases}
			D_t^\varphi v +\nabla^\varphi q= (F_k\cdot \nabla^\varphi) F_k \qquad & \text{in } \Omega,\\
			\nabla^\varphi \cdot v=0 \qquad & \text{in } \Omega,\\
			D_t^\varphi F_j=(F_j\cdot \nabla^\varphi) v  \qquad & \text{in } \Omega,\\
			\nabla^\varphi \cdot  F_j=0 \qquad & \text{in } \Omega\\
			\partial_t \psi=v\cdot N \qquad &\text{on } \Sigma,\\
			q=0 \qquad &\text{on } \Sigma,\\
			F_j\cdot N=0 \qquad &\text{on } \Sigma,\\
			v_3=0 \qquad &\text{on } \Sigma_b,\\
			F_{3j}=0 \qquad &\text{on } \Sigma_b,\\
			q=0 \qquad &\text{on } \Sigma_b,\\
		\end{cases}
	\end{equation}
	with the boundary conditions
	
	\begin{equation}\label{limitbc}
		\begin{cases}
			\partial_t \psi=v\cdot N \qquad &\text{on } \Sigma,\\
			q=0 \qquad &\text{on } \Sigma,\\
			F_j\cdot N=0 \qquad &\text{on } \Sigma,\\
			v_3=0 \qquad &\text{on } \Sigma_b,\\
			F_{3j}=0 \qquad &\text{on } \Sigma_b,\\
			q=0 \qquad &\text{on } \Sigma_b.\\
		\end{cases}
	\end{equation}
	
	To differentiate between the the solution of (\ref{limitsys})-(\ref{limitbc}) and (\ref{phisys})-(\ref{phibc}), we denote the solution of (\ref{phisys})-(\ref{phibc}) by $(\psi^\sigma, v^\sigma, F_j^\sigma)$ to emphasise the dependence of the solution on $\sigma$. We have shown that $|\psi^\sigma(t)|_4^2$ $||v^\sigma(t)||_4^2 , ||F_j^\sigma(t)||_4^2 $ are bounded uniformly in $\sigma$, provided the Rayleigh-Taylor sign condition holds. Hence the Morrey-type embeddings imply that $v^\sigma(t),F_j^\sigma(t)$ are equicontinuous and uniformly bounded in $C^2(\Omega)$ and $\psi^\sigma(t)$ is equicontinuous and uniformly bounded in $C^2(\Sigma)$, which implies $(\psi^\sigma, v^\sigma, F_j^\sigma)\to (\psi, v, F_j)$ in $C^0([0,t],C^2(\Sigma)\times C^2(\Omega)\times C^2(\Omega))$ as $\sigma$ tends to zero.

\bibliographystyle{plain}
\bibliography{Reference.bib}

\end{document}